\documentclass[11pt]{article}
\usepackage{amssymb}
\usepackage{epsfig}
\usepackage{amsmath}
\usepackage{natbib}
\usepackage{latexsym}
\usepackage{amsfonts}
\usepackage{color}
\usepackage{graphicx}

\numberwithin{equation}{section}

\parskip \baselineskip
\newtheorem{lemma}{Lemma}
\newtheorem{proposition}{Proposition}
\newtheorem{theorem}{Theorem}
\newtheorem{corollary}{Corollary}
\newtheorem{remark}{Remark}
\newtheorem{definition}{Definition}

\begin{document}

\begin{center}
{\Large \bf On statistical model extensions based on randomly stopped extremes}
\end{center}

\begin{center}
{\large Jordi Valero, Josep Ginebra\footnote{
       Address for correspondence: Department of Statistics and O.R., Polytechnic University of Catalonia,
       Avgda. Diagonal 647, 6$^{\hbox{\rm a}}$ Planta, 08028Barcelona, Spain
       (E-mail:
       josep.ginebra@upc.es
       )}
        }
\end{center}

The maxima and the minima of a randomly stopped sample of a random variable, $X$,
together with two newly defined random variables
that make $X$ into the maxima or minima of a randomly stopped sample of them, can be used to define
statistical model transformation mechanisms.
These transformations can be used to define models for extreme value data that
are not grounded on large sample theory.
The relationship between the stopping model and
characteristics of the corresponding model transformations obtained is investigated.
In particular, one looks into which stopping models make these model transformations into model extensions,
and which stopping models lead to statistically stable extensions in the sense that
using the model extension a second time leaves the extended model unchanged.
The stopping models under which
the extensions based on randomly stopped maxima and their inverses coincide with
the extensions based on randomly stopped minima and their inverses are also characterized.
The advantages of using models obtained through these model extension mechanisms
instead of resorting to extreme value models grounded on asymptotic arguments
is illustrated by way of examples.



\section{Introduction}		

In disciplines like
hydrology, meteorology, ecology, seismology, actuarial sciences, civil engineering or finance,
there is a need for statistical models to analyze extreme valued data, like the largest one-event rainfall
or the magnitude of the strongest earthquake in a year. In these settings researchers most often resort to
the use of the generalized extreme value model, which is grounded on large sample theory
that only applies as an approximation when sample sizes are large enough.

Hence, there is a need for statistical models for extreme valued data that can be grounded on finite sample theory.
One framework that provides that ground, models the number of events in a year, like the number of rainfalls or of earthquakes,
through a random variable, $N$, with a given stopping model, it models the magnitude of the events in that year as a
sample of i.i.d. observations, $(X_1, \ldots, X_N)$, with a given stopped model,
and it assumes that one observes the maxima or
the minima, $Y$, of that randomly stopped sample.
Models defined like the one for $Y$
are also useful in reliability, where the minimum (or maximum) of a randomly stopped sample
from a lifetime distribution serves as a model for the lifetime of a series (or parallel) system.

Marshall and Olkin (1997) obtains statistical models of this kind by extending an initial
statistical model through the distribution of the
minimum and of the maximum of a geometrically stopped sample of
independent observations with a distribution in the initial family.
This statistical model transformation mechanism
has proved extremely fruitful in practice, as the more than seventeen hundred
citations of that paper indicate.

One nice feature of model transformations based on geometrically stopped extremes is
that they always work as model extensions, because the initial family of distributions is always included in the new family.
A second interesting feature of these geometrically stopped extreme extensions
is that they are statistically stable in the sense that
the extended model can not be further extended by using that same extension mechanism a second time.
These two features are not in place in general,
when transforming statistical models through randomly stopped maxima and minima
with a stopping model different from geometric. In fact, Marshall and Olkin (1997)
conjectures that this kind of stability can only be obtained through geometrically stopped extremes.

Here these issues are investigated in full generality, by looking into all model transformations defined
through the maxima or the minima of $N$-stopped random samples of $X$, for any given stopping model for $N$
and any given stopped model for $X$.

On top of looking into randomly stopped extreme model extensions beyond geometric stopping, we also propose
two new model transformation mechanisms based on two new random variables
defined to be the ones that make $X$ into the randomly stopped maxima and the randomly stopped minima of them, which
we label as the $N$-maxprecursor and the $N$-minprecursor of $X$. These transformations can be viewed
as the inverse transformations of $N$-stopped maxima and of $N$-stopped minima of $X$,
and the statistical models obtained through them can be used
to learn about the magnitude of events, $X$, based on their frequency $N$ and their
extreme values $Y$.

Finally, on top of these four basic model transformation mechanisms based on randomly stopped maxima and minima and on their inverses,
we also propose another two new pair of transformation mechanisms that
combine $N$-stopped maxima of $X$ with their inverses, and combine $N$-stopped minima of $X$ with their inverses.
Under geometric stopping, these combined model transformation mechanisms coincide with the Marshall-Olkin extension mechanism,
and they work as model extensions under any stopping model, which is why we consider them to be the
natural way to generalize Marshall-Olkin when the stopping model is not geometric.

The relationship between characteristics of the stopping model and characteristics
of all the corresponding model transformations considered is studied. The first objective is to identify which
stopping models lead to transformations that always work as model extensions, and the second objective is to identify which
stopping models lead to model extensions that are statistically stable, in the sense that they do not further extend the initial
model beyond the first use.

The second objective leads to the investigation of the class of
stopping models that are closed under pgf composition, because that is a necessary condition for the corresponding
randomly stopped extreme extensions to be stable.
This investigation helps us to disprove by way of examples the
conjecture that only geometric stopping models lead to statistically stable extensions.

The paper also looks into the reversibility conditions
required of stopping models so that the model extensions built based on $N$-stopped maxima and their inverses coincide with the
model extensions built based on $N$-stopped minima and their inverses, which is a property satisfied in particular by
the extensions based on geometric stopping.

The paper illustrates through examples the advantages of modeling extreme valued data with
models obtained through
randomly stopped extreme extensions instead of resorting to the usual generalized extreme
value model backed through large sample arguments. We also use examples to help understand
the rationale behind the use of the models obtained through the new model extension mechanisms that use
the inverse of $N$-stopped maxima or minima.

The paper is organized as follows. Section 2 defines randomly stopped extreme and extreme precursor random variables, and presents
the four basic and four combined model transformation mechanisms that will be investigated, and
Section 3 illustrates the use of models obtained with these transformations to deal with extreme value data.
Section 4 introduces the definition of statistically stable model transformation. Section 5
defines extreme reversible and auto-reversible stopping models, and
Section 6 looks into stopping models that are closed under pgf composition, which are the ones that yield
statistically stable transformations. Section 7 relates these and other
properties of the stopping model with features of the corresponding model transformations,
and Section 8 presents examples of statistically stable randomly stopped extreme extensions.

\section{Statistical models based on randomly stopped extremes}

\subsection{Randomly stopped extremes and extreme precursors}

Let $X$ be a real valued random variable defined through its
cumulative distribution function, $F_X$, and let $N$ be a positive integer valued random variable,
with $\Pr(N=0)=0$, defined through its
probability generating function, $h_N$. Assume that one observes $n$ independent copies of $X$, $X_i$,
where $n$ is a realization of the stopping variable $N$ independent of the $X_i$.

The \emph{$N$-stopped maximum of} $X$, which we denote by $max_N(X)$,
is the random variable $Y = max (X_1,\cdots X_N)$ with cumulative distribution function:
$$ F_{max_N(X)} = h_{N}(F_X), $$
and the \emph{$N$-stopped minimum of} $X$, which we denote by $min_N(X)$,
is the random variable $Y = min (X_1,\cdots, X_N)$ with survival function $h_N(S_X)$, where $S_X=1-F_X$,
and therefore with cdf:
$$  F_{min_N(X)} = 1-h_{N}(1-F_X) = \overline{h}_N(F_X),  $$
where $\overline{h}_N(t)=1 - h_N(1-t)$, which
will be denoted as the \emph{conjugate function} of $h_N(t)$.

These two random variables are studied for example in
Raghundanan and Patil (1972),
Shaked (1975), Consul (1984), Gupta and Gupta (1984),
Rohatgi (1987),
Shaked and Wong (1997),
in pp.155-157 of Arnold et al.(1992) and in
Louzada et al.(2012),

Next, two new random variables that play a central role in what follows are introduced.
They arise from the fact that given any $N$ and any $X$, one can always
interpret $X$ to be the $N$-stopped maximum and the $N$-stopped minimum of the two random variables
defined next.
\begin{definition}
Given any stopping variable $N$ and any real valued random variable $X$ as defined above, let the $N$-maxprecursor of $X$, denoted
as $max_N^{-1}(X)$, be the random variable $Y$ with cdf:
$$  F_{max_N^{-1}(X)} = h_N^{-1}(F_X),  $$
and let the $N$-minprecursor of $X$, denoted as $min_N^{-1}(X)$, be the random variable $Y$ with cdf:
$$ F_{min_N^{-1}(X)}  = \overline{h}_N^{-1}(F_X).  $$
\end{definition}
The properties of $h_N$ and of $\overline{h}_N$ presented in Section 5.1 guarantee that they are always invertible
and therefore that $F_{max_N^{-1}(X)}$ and $F_{min_N^{-1}(X)}$  are always properly defined cdf's.
As a consequence, the random variables $max_N^{-1}(X)$ and $min_N^{-1}(X)$ will exist for any $N$ and any $X$.

By definition, $X$ is always the $N$-stopped maximum of $max_N^{-1}(X)$, the $N$-stopped minimum of $min_N^{-1}(X)$,
the $N$-maxprecursor of $max_N(X)$, and the $N$-minprecursor of $min_N(X)$,
$$  X = max_N(max_N^{-1}(X)) = max_N^{-1}(max_N(X)) = min_N(min_N^{-1}(X)) = min_N^{-1}(min_N(X)),   $$
which is why we denote $N$-maxprecursors and minprecursors as the inverses of the $N$-stopped maxima and minima.

\subsection{Statistical model transformations based on randomly stopped extremes}

Let the family of distributions $\mathcal{X} = \{X_{\theta}: F_{X_{\theta}}, ~ \theta \in \Theta \}$ be a statistical model defined on
$x \in S \subseteq \mathbb{R}$, with parameter space $\Theta$,
where $F_{X_{\theta}}$ is the cdf of $X_{\theta}$.

Let $\mathcal{N}= \{N_{\delta}: h_{N_{\delta}}=\sum_{n=1}^\infty p_n(\delta) t^n, ~ \delta \in \mathcal{D} \}$
be a statistical model
defined on the positive integers, $n \in \mathbb{N}^{+}$, with parameter space $\mathcal{D}$,
where $p_n(\delta) = \Pr(N_{\delta}=n)$ and where $h_{N_{\delta}}$ is the pgf of $N_{\delta}$.
We denote $\mathcal{N}$ as the stopping model.

Note that by definition in this paper it will always be assumed that stopping models, $\mathcal{N}$,
are always such that $\Pr(N_{\delta}=0) = 0$ for any $\delta \in \mathcal{D}$.

We next define four basic mechanisms,$\mathcal{T}(\cdot)$, that transform the initial statistical model, $\mathcal{X}$,
into a new statistical model, $\mathcal{Y}=\mathcal{T}(\mathcal{X})$,
through the $N$-stopped maximum (minimum) of $X \in \mathcal{X}$, and through the $N$-maxprecursors ($N$-minprecursors)
of $X \in \mathcal{X}$, with $N \in \mathcal{N}$.
\begin{definition}
Given any statistical model $\mathcal{X}$ and any stopping model $\mathcal{N}$ as defined above, let
$max_{\mathcal{N}}(\mathcal{X})$ and $max^{-1}_{\mathcal{N}}(\mathcal{X})$ denote the statistical models defined as:
$$ max_{\mathcal{N}}(\mathcal{X}) =
\{Y_{\theta,\delta}: F_{Y_{\theta,\delta}} = h_{N_{\delta}}(F_{X_{\theta}}), ~ \theta \in \Theta, \delta \in \mathcal{D} \}, $$
$$  max^{-1}_{\mathcal{N}}(\mathcal{X}) =
\{Y_{\theta,\delta}: F_{Y_{\theta,\delta}} = h_{N_{\delta}}^{-1}(F_{X_{\theta}}), ~ \theta \in \Theta, \delta \in \mathcal{D} \}.   $$
Likewise, let $min_{\mathcal{N}}(\mathcal{X})$ and $min^{-1}_{\mathcal{N}}(\mathcal{X})$ denote the statistical models defined as:
$$   min_{\mathcal{N}}(\mathcal{X}) =
\{Y_{\theta,\delta}: F_{Y_{\theta,\delta}}=\overline{h}_{N_{\delta}}(F_{X_{\theta}}), ~ \theta \in \Theta, \delta \in \mathcal{D} \},  $$
$$   min^{-1}_{\mathcal{N}}(\mathcal{X}) =
\{Y_{\theta,\delta}: F_{Y_{\theta,\delta}}=\overline{h}_{N_{\delta}}^{-1}(F_{X_{\theta}}), ~ \theta \in \Theta, \delta \in \mathcal{D} \}. $$
\end{definition}
These two pairs of basic transformations do not always work as model extensions.
Instead, the two pairs of combined transformations defined next work as model extensions for any $\mathcal{X}$,
even when one of the two new parameters is fixed. They are the family of
all $N$-stopped maximum (minimum) of all $N$-maxprecursors ($N$-minprecursors) of $X$, and viceversa.
\begin{definition}
Given any statistical model $\mathcal{X}$ and any stopping model $\mathcal{N}$ as defined above, let
$max_{\mathcal{N}}(max^{-1}_{\mathcal{N}}(\mathcal{X}))$ and $max^{-1}_{\mathcal{N}}(max_{\mathcal{N}}(\mathcal{X}))$
denote the statistical models defined as:
$$  max_{\mathcal{N}}(max^{-1}_{\mathcal{N}}(\mathcal{X})) = \{Y_{\theta,\delta_1,\delta_2}: F_{Y_{\theta,\delta_1,\delta_2}} =
h_{N_{\delta_2}}\circ h_{N_{\delta_1}}^{-1}(F_{X_{\theta}}), ~ \theta \in \Theta, \delta_1,\delta_2 \in \mathcal{D} \},  $$
$$  max^{-1}_{\mathcal{N}}(max_{\mathcal{N}}(\mathcal{X})) = \{Y_{\theta,\delta_1,\delta_2}: F_{Y_{\theta,\delta_1,\delta_2}} =
h_{N_{\delta_2}}^{-1}\circ h_{N_{\delta_1}}(F_{X_{\theta}}), ~ \theta \in \Theta, \delta_1,\delta_2 \in \mathcal{D} \}. $$
Likewise, let $min_{\mathcal{N}}(min^{-1}_{\mathcal{N}}(\mathcal{X}))$ and $min^{-1}_{\mathcal{N}}(min_{\mathcal{N}}(\mathcal{X}))$
denote the statistical models:
$$  min_{\mathcal{N}}(min^{-1}_{\mathcal{N}}(\mathcal{X})) = \{Y_{\theta,\delta_1,\delta_2}: F_{Y_{\theta,\delta_1,\delta_2}} =
\overline{h}_{N_{\delta_2}}\circ \overline{h}_{N_{\delta_1}}^{-1}(F_{X_{\theta}}), ~ \theta \in \Theta, \delta_1,\delta_2 \in \mathcal{D} \}, $$
$$  min^{-1}_{\mathcal{N}}(min_{\mathcal{N}}(\mathcal{X})) = \{Y_{\theta,\delta_1,\delta_2}: F_{Y_{\theta,\delta_1,\delta_2}} =
\overline{h}_{N_{\delta_2}}^{-1}\circ \overline{h}_{N_{\delta_1}}(F_{X_{\theta}}), ~ \theta \in \Theta, \delta_1,\delta_2 \in \mathcal{D} \}. $$
\end{definition}
Note that these statistical model transformation mechanisms can also be used to generate statistical models, $\mathcal{Y}$,
starting from a single initial random variable, $\mathcal{Y} = \mathcal{T}(X)$.

Using our notation, the Marshall-Olkin extension of $\mathcal{X}$
is defined to be $max_{\mathcal{N}}(\mathcal{X}) \cup min_{\mathcal{N}}(\mathcal{X})$ when $\mathcal{N}$ is
the geometric stopping model. In Sections 7.2 and 7.4 it will be argued that for stopping models other than geometric
the transformation $max_{\mathcal{N}}(\mathcal{X}) \cup min_{\mathcal{N}}(\mathcal{X})$
does not always work as an extension,
but that under geometric stopping this transformation coincides with the four model transformations
in Definition 3, which do work as extensions under any stopping model.
As a consequence, we will propose Definition 3 and not $max_{\mathcal{N}}(\cdot) \cup min_{\mathcal{N}}(\cdot)$
to be the natural way to generalize Marshall-Olkin when using stopping models different from geometric.

\section{Examples of the use of randomly stopped extreme models}

The examples presented here illustrate the advantage in using
models defined through the randomly stopped extreme transformations in Definition 2
instead of using the generalized extreme value model, and they help understand
the practical relevance of the randomly stopped extreme precursor models also considered in that definition.
The examples also touch on the rationale behind the use of the model extensions
proposed in Definition 3.

\subsection{On the usefulness of randomly stopped extreme models}

Let's assume for example that one has data on the rainfall in the largest rain event of a year, $Y_i$,
for a set of $m$ years, $(y_1,\ldots,y_m)$. This kind of data is usually modeled through
the three parameter generalized extreme value model,
because it is the limiting model for properly normalized extreme values
when the rainfall in an event is i.i.d., and the number of rainfall events in a year grows.

As an alternative way to model this kind of data one can assume that
the number of rain events in the $i$-th year, $N_i$, is random and can be modeled through a specific stopping model,
$\mathcal{N}$, and that the rainfall in the set of $N_i$ events is a sample of i.i.d. observations, $(X_1,\ldots,X_{N_i})$,
from a specific model, $\mathcal{X}$.
In this framework the statistical model for the largest rainfall in the $i$-th year,
$Y_i = \max{(X_1,\ldots,X_{N_i})}$, is
the $\mathcal{Y} = max_{\mathcal{N}}(\mathcal{X})$
considered in Definition 2 for that $\mathcal{N}$ and that $\mathcal{X}$.

In particular, for simplicity here it will be assumed that the stopping model for the number of rain events, $N_p$,
is the Logarithmic$(p)$ model covered in Example 5.2 and in Appendix 1,
and that the model for the rainfall in an event, $X_{\lambda}$, is the Exponential$(\lambda)$ with cdf
$F_{X_{\lambda}}(x)=1-e^{-\lambda x}$.
In that case, the model for the largest rainfall of the year, $(Y_1, \ldots, Y_m)$,
is the logarithmic stopped maximum of an exponential,
$$  \mathcal{Y}_{Lg-Exp} =
\{Y_{p,\lambda}: F_{Y_{p,\lambda}} = h_{N_p}(F_{X_{\lambda}}) = \frac{\log{(1-p+p e^{-\lambda y})}}{\log{(1-p)}}    , \lambda \in (0,\infty), p \in (0,1) \}.    $$
To compare the use of the randomly stopped extreme models with the use of the generalized extreme value model,
we have simulated a sample for $m=150$ years
assuming that $N_i$ is Logarithmic$(p=.95)$ and $X_i$ is Exponential$(\lambda=0.01)$.
We have fitted the true two-parameter $\mathcal{Y}_{Lg-Exp}$ model and
the three parameter generalized extreme value model,
$$  \mathcal{Y}_{GEV} = \{Y_{\eta,\theta,\kappa}: F_{Y_{\eta,\theta,\kappa}} =
e^{-(1-\kappa(x-\eta)/\theta)^{1/\kappa}}, \eta \in (-\infty,\infty) , \theta \in (0,\infty), \kappa \in (-\infty,\infty) \},  $$
on this data set by maximum likelihood. We have also fitted the
$\mathcal{Y}_{TB2-Exp}$ and $\mathcal{Y}_{ETNB-Exp}$ models, which are
the randomly stopped maxima of an exponential sample
when the stopping model is the truncated binomial$(2,p)$ and the extended truncated negative binomial (ETNB)
with pgf $h_N = \frac{\log{(1-pt)^{-r}-1}}{\log{(1-p)^{-r}-1}}$ where $r$ is in $(-1,\infty)$.
We also fit $\mathcal{Y}_{PC-LgNor}$, which is the randomly stopped maximum with $\mathcal{N}$
being the potential conjugate model considered in Example 6.1
and $\mathcal{X}$ being the lognormal model.

Table 1 presents the maximum likelihood estimates of the parameters of these five
models together with the value of the log-likelihood at its maximum, and their AIC and BIC.
Note that the $\mathcal{Y}_{ETNB-Exp}$ model fits the data slightly better than
the actual $\mathcal{Y}_{Lg-Exp}$ model, but when $r=0$ the $\mathcal{Y}_{ETNB-Exp}$ becomes the $\mathcal{Y}_{Lg-Exp}$ model and
the likelihood ratio test between these two nested models
does not reject the simpler actual model with a $p-val$ of $0.758$.

\begin{table}[tbp]\centering
{\small
\begin{tabular}{|c|c|ccc|c|c|c|}
\cline{1-8}
Model & N.par & & MLE &  & loglikel  & AIC & BIC \\
\cline{1-8}
$\mathcal{Y}_{Lg-Exp}$ & $2$ & $\hat{p}=.9574$ &  & $\hat{\lambda}=.0098$ & $-942.326$ & $1888.65$ & $1894.67$  \\
\cline{1-8}
$\mathcal{Y}_{ETNB-Exp}$ & $3$ & $\hat{p}=.9844$ & $\hat{r}=-.1177$   & $\hat{\lambda}=.0108$ &  $-942.172$ & $1890.34$ & $1899.38$ \\
\cline{1-8}
$\mathcal{Y}_{TB2-Exp}$ &  $2$ & $\hat{p}=.3949$  &  & $\hat{\lambda}=.0063$ & $-945.196$  & $1894.39$  & $1900.41$  \\
\cline{1-8}
$\mathcal{Y}_{PC-LgNor}$ & $3$ & $\hat{p}=.9352$ & $\hat{\mu}=4.9109$  & $\hat{\sigma}=1.1475$ & $-952.578$ & $1911.15$ & $1920.19$ \\
\cline{1-8}
$\mathcal{Y}_{GEV}$ &  $3$ & $\hat{\eta}=129.01$ & $\hat{\theta}=111.25$ & $\hat{\kappa}=-.1207$ & $-954.256$  & $1914.51$ & $1923.54$ \\
\cline{1-8}
\end{tabular}
 \caption{Maximum likelihood parameter estimates, logarithm of the likelihood at the mle,
and AIC and BIC for the five models considered for the data on the largest annual rainfall event.}}
 \label{t2}
\end{table}

Even though the $\mathcal{Y}_{GEV}$ model has one more parameter than the
actual $\mathcal{Y}_{Lg-Exp}$ model, it fits the simulated
data significantly worse than this model, and worse than the other three stopped extreme models tried,
even though two of these models assume a wrong stopping model and one of them assumes a wrong
stopping and a wrong stopped model.
Of course that will not always be the case,
and the $\mathcal{Y}_{GEV}$ model will do better than other stopped extreme models, but when one
has a good guess on what the stopping and the stopped models could be, the corresponding
randomly stopped extreme model will tend to do better than $\mathcal{Y}_{GEV}$.

Note also that an important advantage of using randomly stopped extreme models is that through them
one can interpret the estimated parameter values in terms of the parameters
of the model for the stopping variable and the parameters of the model for the stopped variable.
That provides useful information about the frequency of rain and
about the distribution of the amounts of rain in them,
which is lacking when the analysis is based on the GEV model.

\textbf{Remark:} When $\mathcal{Y} = max_{\mathcal{N}}(\mathcal{X})$ one has that
$F_Y=h_N(F_X)$ and the pdf of $Y$ is $f_Y=h_N^{'}(F_X) f_X$, where $f_X$ is the pdf of $X$.
Therefore $f_Y$ is a weighted version of $f_X$ and maximizing
the likelihood function using data on $Y$ is not any more complicated than doing it with data on $X$.

\subsection{On the usefulness of randomly stopped extreme precursors}

Lets assume here that one has data on the magnitude of the strongest earthquake on a given year
for $m_1$ years, $(y_1,\ldots,y_{m_1})$, and data on the number of earthquakes
in a year for $m_2$ years, $(n_1, \ldots, n_{m_2})$, where the set of years with available data might not coincide.
Let's also assume that one has a good model $\mathcal{Y}$ for
$Y_i$ and a good model $\mathcal{N}$ for $N_i$.

Like in the previous example one can pose
$Y_i = \max{(X_1,\ldots,X_{N_i})}$ where the magnitudes of the earthquakes, $X_j$, are i.i.d.
realizations of a random variable, $X$, and hence one can assume that
$\mathcal{Y}=max_{\mathcal{N}}(\mathcal{X})$. In such a setting
one might lack data about the $X_j$ and yet the interest in the analysis might be to learn about the
distribution of these $X_j$, and therefore  about their cdf, $F_{X}$.

In particular, the stopping model for the number of earthquakes, $N_i$, could for example again
be Logarithmic$(p)$, and a good model for the magnitude of the strongest earthquake, $Y_i$,
could be the GEV$(\eta,\theta,\kappa)$ model that was discarded in the previous example
for the largest rainfall.
If that was the case the magnitude of earthquakes, $X_j$, would be a sample from
the random variable $X$ that is the $N_p$-maxprecursor of the GEV r.v., $Y_{\eta,\theta,\kappa}$, and
the cdf of $X$ would be:
$$ F_{X_{p,\eta,\theta,\kappa}} = F_{max_{N_p}^{-1}(Y_{\eta,\theta,\kappa})} = h_{N_p}^{-1}(F_{Y_{\eta,\theta,\kappa}}).  $$
Hence, by obtaining maximum likelihood estimates of $p$ and of $(\eta,\theta,\kappa)$ and
estimates of their standard deviations using the data on
$N_i$ and the data on $Y_i$ one would obtain estimates and confidence intervals for the cdf of $X$,
$\hat{F}_{X_{p,\eta,\theta,\kappa}} = h_{N_{\hat{p}}}^{-1}(F_{Y_{\hat{\eta},\hat{\theta},\hat{\kappa}}})$.

\subsection{On the rationale behind using the extensions in Definition 3}

Finally, lets assume that in either the hydrology or the seismology settings considered above one guesses that
$\mathcal{N}_0$ is the stopping model for the number of events, $N_i$, and
$\mathcal{X}_0$ is the model for the magnitude of the events $X_j$,
but it turns that the statistical model $\mathcal{Y}_0 = max_{\mathcal{N}_0}(\mathcal{X}_0)$ for $Y_i = \max{(X_1,\ldots,X_{N_i})}$
fails to fit properly the sample of extreme values available, $(y_1,\ldots,y_m)$.

In a case like this, if one is confident that $\mathcal{N}_0$ is the right stopping model
one will want to extend $\mathcal{Y}_0$ by extending $\mathcal{X}_0$
while still using $\mathcal{N}_0$ as the stopping model.
The first model extension in Definition 3 does that by replacing
$\mathcal{Y}_0 = max_{\mathcal{N}_0}(\mathcal{X}_0)$ by:
$$  \mathcal{Y}_1 = max_{\mathcal{N}_0}(max^{-1}_{\mathcal{N}_0}(\mathcal{Y}_0)) =
\{Y_{\xi,\delta_1,\delta_2}: F_{Y_{\xi,\delta_1,\delta_2}} =
h_{N_{\delta_2}}\circ h_{N_{\delta_1}}^{-1}(F_{Y_{\xi}}), ~ \xi \in \Xi, \delta_1,\delta_2 \in \mathcal{D} \},  $$
where $\Xi$ is the parameter space of $\mathcal{Y}_0$.
In this way, the extended model can be posed as $\mathcal{Y}_1 = max_{\mathcal{N}_0}(\mathcal{X}_1)$
where $\mathcal{X}_0$ has been replaced by its extension,
$\mathcal{X}_1=max^{-1}_{\mathcal{N}_0}(max_{\mathcal{N}_0}(\mathcal{X}_0))$.
Note that this extension also applies when $\mathcal{Y}_0$ is chosen
without making any $\mathcal{X}_0$ explicit,
in which case the extended model is $\mathcal{Y}_1 = max_{\mathcal{N}_0}(\mathcal{X}_1)$
with $\mathcal{X}_1=max^{-1}_{\mathcal{N}_0}(\mathcal{Y}_0)$.

By construction, the dimension of the parameter space of
$\mathcal{Y}_1=max_{\mathcal{N}_0}(max^{-1}_{\mathcal{N}_0}(\mathcal{Y}_0))$ is never smaller than
the one of $\mathcal{Y}_1^{'}=max^{-1}_{\mathcal{N}_0}(\mathcal{Y}_0)$, which is never smaller than
the one of $\mathcal{Y}_0$. This paper investigates when is the initial model always included in the
transformed model, and when does repeated use of these extensions fail to keep extending the model.

\section{Statistical stability of statistical model transformations}

Transformations of a statistical model, $\mathcal{X}$, into
a new model, $\mathcal{Y}=\mathcal{T}(\mathcal{X})$, can be classified depending on how
initial and final models relate. Most often neither $\mathcal{X}$ nor $\mathcal{Y}$ are included into
each other. The next definition distinguishes three
possible relationships when they do.
\begin{definition}
Let $\mathcal{T}(\cdot)$ transform a statistical model, $\mathcal{X}$, into
$\mathcal{Y}=\mathcal{T}(\mathcal{X})$. Then
\begin{enumerate}\vspace{-.3cm}
\item{}
if $\mathcal{X} \subset  \mathcal{T}(\mathcal{X})$,
one says that $\mathcal{X}$ is extended by $\mathcal{T}(\cdot)$, and that $\mathcal{T}(\cdot)$ extends $\mathcal{X}$,
\item{}
if $\mathcal{T}(\mathcal{X}) \subset \mathcal{X}$,
one says that $\mathcal{X}$ is contracted by $\mathcal{T}(\cdot)$, and that $\mathcal{T}(\cdot)$ contracts $\mathcal{X}$,
\item{}
if $\mathcal{T}(\mathcal{X}) = \mathcal{X}$,
one says that $\mathcal{X}$ is invariant under $\mathcal{T}(\cdot)$.
\end{enumerate}\vspace{-.3cm}
\end{definition}
When $\mathcal{X}$ is extended by $\mathcal{T}(\cdot)$ for all $\mathcal{X}$, one says that $\mathcal{T}(\cdot)$
is a model extension. Most often, using a model extension repeatedly will keep extending the model, but some model
extensions do not further extend models beyond their first use.
These special model extensions are examples of the statistically stable transformations defined next.
\begin{definition}
A statistical model transformation, $\mathcal{T}(\cdot)$, is said to be statistically stable
if for any model $\mathcal{X}$ one has that $\mathcal{T}(\mathcal{X})$
is invariant under $\mathcal{T}(\cdot)$, and so if $\mathcal{T}(\mathcal{T}(\mathcal{X}))=\mathcal{T}(\mathcal{X})$ for any $\mathcal{X}$.
\end{definition}
When a model transformation is statistically stable, using that transformation twice in a row
on any statistical model, $\mathcal{X}$, has the same effect as using it just once.

Definition 5 generalizes to any statistical model transformation the concept of
geometric-extreme stability proposed in Marshall and Olkin (1997) in the special case of geometric stopped extreme transformations.
Note that the statistical notion of stability presented here
is different from probabilistic notions of stability, like the ones used in
Rachev and Resnick (1991) or in Fama and Roll (1968), which apply to individual random variables and not
to families of them.

The main purpose of the paper is to investigate the properties of the model transformations in Definitions 2 and 3, and to
determine when do they work as model extensions, and when are these model extensions statistically stable in the sense of Definition 5.
This depends only on the characteristics of the stopping model, and in particular on whether they are
extreme auto-reversible and/or closed under pgf composition, the way defined in the next two sections.

\section{Stopping models that are extreme reversible or auto-reversible}

\subsection{Properties of $h_N$, $\overline{h}_N$, $h_N^{-1}$, and $\overline{h}_N^{-1}$ for positive count variables}

A function, $h_N$, is the probability generating function of a positive integer-valued random variable $N$,
if and only if it is real valued and such that
$h_N(0)=0$, that $h_N(1)=1$, and that
it is analytic at least on $[0,1)$,
with all derivatives in that set being non-negative.

As a consequence, $\overline{h}_N(t)=1 - h_N(1-t)$
is always such that $\overline{h}_N(0)=0$, $\overline{h}_N(1)=1$, and that
it is analytic at least on $(0,1]$, with all of its odd derivatives in that set being
non-negative, and all of its even derivatives non-positive.
If all the moments of $N$ are
finite, analyticity and the declared signs of the derivatives of $h_N$ and $\overline{h}_N$ hold at least on $[0,1]$.

From the characterization of $h_N$ it also follows that $h_N^{-1}$ and $\overline{h}_N^{-1}$ are always such
that $h_N^{-1}(0) = \overline{h}_N^{-1}(0) = 0$ and $h_N^{-1}(1) = \overline{h}_N^{-1}(1) = 1$, and
they are analytic at least on $(0,1)$ with a first derivative that is non-negative in that set.
The second derivative of $h_N^{-1}$ is non-positive, while the second derivative of
$\overline{h}_N^{-1}$ is non-negative.

In particular, $h_N$, $\overline{h}_N$, $h_N^{-1}$ and $\overline{h}_N^{-1}$
are always continuous and increasing on $[0,1]$, with
$h_N$ and $\overline{h}_N^{-1}$ being convex,
and $\overline{h}_N$ and $h_N^{-1}$ being concave.

For the limiting stopping random variable
$N_I$ with $\Pr(N_I=1)=1$, these four functions coincide,
$h_{N_I}(t) = t = \overline{h}_{N_I}(t) = h_{N_I}^{-1}(t) = \overline{h}_{N_I}^{-1}(t)$.
The next result will be used later on.
\begin{proposition}
If $N, N_1$ and $N_2$ are positive integer valued random variables with pgfs $h_N$, $h_{N_1}$ and $h_{N_2}$, then
1) $\overline{\overline{h}}_N = h_N,$
2) $\overline{h_N^{-1}} = \overline{h}_N^{-1},$ and
3) $\overline{h}_{N_1}\circ\overline{h}_{N_2} = \overline{(h_{N_1}\circ h_{N_2})}.$
\end{proposition}

\subsection{Extreme reversible stopping models}

As a consequence of the properties listed above, $\overline{h}_N$ and $h_N^{-1}$ can only be the pgf
of a positive integer valued random variable if $N=N_I$, with $h_{N_I} = t$.

On the other hand,
$\overline{h}_N^{-1}$ sometimes is the pgf of a non-degenerate positive integered random variable, $N^*$.
That leads to the following definition.
\begin{definition}
The pair of positive integer valued random variables, $(N,N^*)$, is said to be extreme reversible
if $\overline{h}_N^{-1} = h_{N^*}$, and therefore if $h_N^{-1} = \overline{h}_{N^*}$.
\end{definition}
When $(N,N^*)$ are extreme reversible, their pgf's need to be such that:
$$ h_{N^*}\circ\overline{h}_{N}(t) =  \overline{h}_{N^*}\circ h_{N}(t) = t = \overline{h}_{N}\circ h_{N^*}(t) = h_{N}\circ \overline{h}_{N^*}(t),
~~ \mbox{for} ~~ t \in [0,1], $$
and in that case, $max_N^{-1}(X) = min_{N^*}(X)$,
$min_N^{-1}(X) = max_{N^*}(X)$, and therefore
$$ X = max_N(min_{N^*}(X)) = min_N(max_{N^*}(X)) = max_{N^*}(min_{N}(X)) = min_{N^*}(max_N(X)). $$
It is important to emphasize that extreme reversibility
is a property of $(N, N^*)$, and that when it holds, this property applies for any real valued random variable, $X$.

\textbf{Example 5.1:}
For the ``potential conjugate" random variable $N_b$, with $h_{N_b} = 1 - (1 - t)^b$
for $b \in (0,1]$, one has that $\overline{h}_{N_b}^{-1}= t^{1/b}$,
which is a pgf when $b=1/m$ and $m$ is a positive integer.
Hence, for any positive integer $m$, the $N_m$ with pgf $h_{N_m} = 1 - (1 - t)^{1/m}$,
and $N_m^*$ with pgf $h_{N_m^*}=t^m$, are extreme reversible.

\textbf{Example 5.2:}
If $N_{\alpha}$ is zero-truncated Poisson($\alpha$), with:
$$ h_{N_{\alpha}} = \frac{e^{\alpha t} - 1}{e^{\alpha} - 1} $$
for a given given $\alpha > 0$, then:
$$ \overline{h}_{N_{\alpha}}^{-1} = - \frac{1}{\alpha} \ln(1 - (1 - e^{-\alpha})t) = h_{N_{\alpha}^{*}}, $$
which is the pgf of a r.v. $N_{\alpha}^{*}$ with a Logarithmic($\alpha$) distribution, most often
parametrized through $p=1-e^{-\alpha}$. This means that each zero-truncated Poisson random variable is extreme reversible
with one logarithmic random variable.

If a statistical model, $\mathcal{N}^*$, is the set of all random variables $N^*$ that are extreme reversible with
a random variable in $\mathcal{N}$, one says that $\mathcal{N}^*$ and $\mathcal{N}$ are a pair of extreme reversible models.

Note that when $\mathcal{N}$ and $\mathcal{N}^*$ are extreme reversible one has that $max_{\mathcal{N}}^{-1}(\cdot) = min_{\mathcal{N}^*}(\cdot)$, and that $min_{\mathcal{N}}^{-1}(\cdot) = max_{\mathcal{N}^*}(\cdot)$, and one also has that:
$$ max_{\mathcal{N}}(max^{-1}_{\mathcal{N}}(\cdot)) = min^{-1}_{\mathcal{N}^{*}}(min_{\mathcal{N}^{*}}(\cdot)), $$
$$ max^{-1}_{\mathcal{N}}(max_{\mathcal{N}}(\cdot)) = min_{\mathcal{N}^{*}}(min^{-1}_{\mathcal{N}^{*}}(\cdot)), $$
and viceversa. As an immediate consequence, when $\mathcal{N}$ and $\mathcal{N}^*$ are extreme reversible
models the set of transformations in Definitions 2 and 3 obtained with $\mathcal{N}$ and the set of transformations in these definitions
obtained with $\mathcal{N}^*$ coincide.

\subsection{Extreme auto-reversible stopping models}

There are instances when $N$ and $N^*$ are the same, hence the next definition.
\begin{definition}
The positive integer random variable $N$ is extreme auto-reversible
if $\overline{h}_N^{-1} = h_{N}$, and therefore if $h_N^{-1} = \overline{h}_{N}$.
\end{definition}
When $N$ is extreme auto-reversible,
$$   h_N\circ \overline{h}_N(t) = \overline{h}_N\circ h_{N}(t)
= t, ~~ \mbox{for} ~~ t \in [0,1], $$
which is a condition used in stochastic comparison theorems of Shaked (1975) and Shaked and Wong (1997). When it holds,
$max_N^{-1}(X) = min_{N}(X)$,
$min_N^{-1}(X) = max_{N}(X)$, and:
$$  X = max_N(min_{N}(X)) = min_N(max_{N}(X)).  $$
A necessary condition for a r.v. $N$ to be auto-reversible is that $\Pr(N=1)=1/E[N]$.
The next result, providing a way to generate two auto-reversible random variables
starting from any pair of reversible ones, will be used to find examples of auto-reversible variables.
\begin{proposition}
If $(N,N^*)$ are a pair of extreme reversible random variables, with pgf's $h_N$ and $h_{N^*} = \overline{h}_N^{-1}$,
then the random variables $N_1$ and $N_2$, with pgfs
$h_{N_1} =  h_N\circ h_{N^{*}}$ and $h_{N_2} = h_{N^{*}}\circ h_N$,
are both extreme auto-reversible.
\end{proposition}
{\sl{Proof:}} Given that $\overline{h}_N(t) = 1 - h_N(1-t)$, one has that:
$$ h_{N_1}\circ \overline{h}_{N_1}(t) =
h_N\circ h_{N^*}\circ \overline{h}_N\circ \overline{h}_{N^*}(t) =  h_N\circ \overline{h}_{N^*}(t) = t,   $$
where the last two steps use the fact that $N$ and $N^*$ are extreme reversible.
\hfill{$\square$}

\begin{corollary} If $N$ is extreme auto-reversible with pgf $h_N$, then the random variable $N_3$ with pgf
$h_{N_3} = h_N\circ h_N$ is also extreme auto-reversible.
\end{corollary}

\textbf{Example 5.3:}
Using Proposition 2 with the random variables of Example 5.1 leads to $h_{N_1} = 1 - (1-t^m)^{1/m}$
and to $h_{N_2} = (1 - (1-t)^{1/m})^m$,
which whenever $m$ is a positive integer are the pgf's of two auto-reversible random variables.

\textbf{Example 5.4:}
Using Proposition 2 with the random variables of Example 5,2 yields:
$$ h_{N_1} = \frac{pt}{1-(1-p)t}, $$
for $0 < p = e^{-\alpha} \le 1$, which is the pgf of the geometric distribution, and
$$ h_{N_2} = 1 - (1/\alpha) \log(1+e^{\alpha} - e^{\alpha t}),  $$
for $\alpha > 0$, where $N_1$ and $N_2$ are extreme auto-reversible random variables.

When all random variables $N$ in $\mathcal{N}$
are extreme auto-reversible, one says that the stopping model $\mathcal{N}$ is extreme auto-reversible.

When $\mathcal{N}$ is an extreme auto-reversible model one has that $max_{\mathcal{N}}^{-1}(\cdot) = min_{\mathcal{N}}(\cdot)$ and that
$min_{\mathcal{N}}^{-1}(\cdot) = max_{\mathcal{N}}(\cdot)$, and therefore that:
$$max^{-1}_{\mathcal{N}}(max_{\mathcal{N}}(\cdot)) = min_{\mathcal{N}}(min^{-1}_{\mathcal{N}}(\cdot)) = min_{\mathcal{N}}(max_{\mathcal{N}}(\cdot)),$$
$$ min^{-1}_{\mathcal{N}}(min_{\mathcal{N}}(\cdot)) = max_{\mathcal{N}}(max^{-1}_{\mathcal{N}}(\cdot)) = max_{\mathcal{N}}(min_{\mathcal{N}}(\cdot)).$$
Therefore, when $\mathcal{N}$ is an extreme auto-reversible model, the four basic and four combined transformations in Definitions 2 and 3
collapse down into two basic and two combined transformations.

\section{Stopping models closed under pgf composition}

A necessary condition for the transformations in Definitions 2 and 3 to be statistically stable is that the
corresponding stopping model are closed under pgf composition as defined next.
\begin{definition}
The stopping model $\mathcal{N}  = \{N_{\delta}: h_{N_{\delta}}, ~ \delta \in \mathcal{D} \}$
is said to be closed under pgf composition, if having
$N_{\delta_1}$ and $N_{\delta_2}$ with pgfs $h_{N_{\delta_1}}$ and $h_{N_{\delta_2}}$ belonging to $\mathcal{N}$, implies that $N_{\delta_3}$
with pgf $h_{N_{\delta_3}} = h_{N_{\delta_1}}\circ h_{N_{\delta_2}}$
also belongs to $\mathcal{N}$.
\end{definition}
Requiring that $\mathcal{N}$ be closed under pgf composition is equivalent to requiring that
if $N_{\delta_1}$ and $N_{\delta_2}$ belong to $\mathcal{N}$, then the $N_{\delta_1}$-stopped sum of $N_{\delta_2}$
also belongs to $\mathcal{N}$, and it is thus equivalent to being closed under model compounding.

\subsection{Uniparametric stopping models closed under pgf composition}

Here we restrict consideration to stopping models,
$ \mathcal{N} = \{N_{\delta}: h_{N_{\delta}}=\sum_{i=1}^\infty p_i(\delta) t^i, ~ \delta \in \mathcal{D} \},$ that
i) are closed under pgf composition,
ii) have a parametrization $\delta$ such that the $p_i(\delta) = \Pr(N_{\delta}=i)$ are
continuously differentiable in $\delta$ for any $i$, and
iii) have a parameter space, $\mathcal{D}$, that is a connected subset of $\mathbb{R}$ with a non-empty interior.
From now on, this class of stopping models
is denoted in a shorthanded way just as ``\emph{models uniparametric and closed under pgf composition}."

By focusing on stopping models continuously differentiable and with this kind of parameter space,
we restrict consideration to the kind of stopping models useful in statistical practice.
In particular, we essentially require that the parameter space be a non-empty interval, thus
avoiding stopping models closed under pgf composition like $\mathcal{N}  = \{N_k: h_{N_k} = t^k, ~ k \in \mathbb{N}^{+} \}$,
which lead to trivially stable transformations, and we also avoid parameter spaces with isolated points.

The following result, crucial in all that follows, is proved in Appendix 2.
\begin{theorem}
If a stopping model, $\mathcal{N}  = \{N_{\delta}: ~ \delta \in \mathcal{D} \}$, is
``\emph{uniparametric and closed under pgf composition}" as defined above,
then:
\begin{enumerate}\vspace{-.3cm}
\item{} $p_1(\delta)=\Pr(N_{\delta}=1) > 0$ for all $N_{\delta} \in \mathcal{N}$,
\item{} $\mathcal{N}$ can be parametrized in an identifiable way through $\theta=\Pr(N_{\delta}=1)$,
or through $\eta = - \log{\Pr(N_{\delta}=1)}$, and
\item{} the parameter space is of the form $(0,\theta_0]$ for a given $\theta_0 \le 1$ when using $\theta$,
and it is of the form $\mathcal{H} = [\eta_0,\infty)$ for a given $\eta_0 \ge 0$ when using $\eta$.
\end{enumerate}\vspace{-.3cm}
\end{theorem}
From now on, we will always use $\eta$ as the parametrization for models \emph{uniparametric and closed under pgf composition}.
Note that $N_I$, with $h_{N_I}(t)=t$, belongs to one of these models if, and
only if, the lower limit of the parameter space, $\eta_0$, is equal to $0$.

Next consequence of Theorem 1 relates to repeated use of the transformations in Definition 2.
\begin{theorem}
If the stopping model, $\mathcal{N} = \{N_{\eta}: h_{N_{\eta}}, ~ \eta \in [\eta_0,\infty) \}$, is
``\emph{uniparametric and closed under pgf composition}" as defined above, then:
\begin{enumerate}\vspace{-.3cm}
\item{} $h_{N_{\eta_1}}\circ h_{N_{\eta_2}} = h_{N_{\eta_2}}\circ h_{N_{\eta_1}} = h_{N_{\eta_1+\eta_2}}$,
\item{} $\overline{h}_{N_{\eta_1}}\circ \overline{h}_{N_{\eta_2}} = \overline{h}_{N_{\eta_2}}\circ \overline{h}_{N_{\eta_1}}
= \overline{h}_{N_{\eta_1+\eta_2}}$,
\item{} $h_{N_{\eta_1}}^{-1}\circ h_{N_{\eta_2}}^{-1} = h_{N_{\eta_2}}^{-1}\circ h_{N_{\eta_1}}^{-1}
= h_{N_{\eta_1+\eta_2}}^{-1}$,
\item{} $\overline{h^{-1}}_{N_{\eta_1}}\circ \overline{h^{-1}}_{N_{\eta_2}}
= \overline{h^{-1}}_{N_{\eta_2}}\circ \overline{h^{-1}}_{N_{\eta_1}}
= \overline{h^{-1}}_{N_{\eta_1+\eta_2}}$,
\end{enumerate}\vspace{-.3cm}
for all $\eta_1, \eta_2 \in [\eta_0,\infty)$.
\end{theorem}
{\sl{Proof:}} Given that $\mathcal{N}$ is
closed under pgf composition, $h_{N_{\eta_1}}\circ h_{N_{\eta_2}} = h_{N_{\eta}}$, with:
$$ \eta = - \log((h_{N_{\eta_1}}(h_{N_{\eta_2}}(t)))_{|t=0}^{\prime}) =
-\log((h_{N_{\eta_1}}^{\prime}(h_{N_{\eta_2}}(t))h_{N_{\eta_2}}^{\prime}(t))_{|t=0}) =
\eta_1 + \eta_2,    $$
and commutativity follows from the commutativity of addition.
The other three assertions follow from the fact that, because of Proposition 1,
$$\overline{h}_{N_{\eta_1}}\circ\overline{h}_{N_{\eta_2}}=
\overline{(h_{N_{\eta_1}}\circ h_{N_{\eta_2}})}=\overline{h}_{N_{\eta_1+\eta_2}}, $$
$$h_{N_{\eta_1}}^{-1}\circ h_{N_{\eta_2}}^{-1}=(h_{N_{\eta_2}}\circ h_{N_{\eta_1}})^{-1}=h_{N_{\eta_1 + \eta_2}}^{-1}, $$
and
$$\overline{h^{-1}}_{N_{\eta_1}}\circ\overline{h^{-1}}_{N_{\eta_2}}=\overline{(h_{N_{\eta_2}}\circ h_{N_{\eta_1}})^{-1}}=\overline{h^{-1}}_{N_{\eta_1 + \eta_2}}. $$
\hfill{$\square$}

The second result that follows from Theorem 1
will imply that under stopping models closed under pgf composition, Definition 3 yields only
two distinct extensions, and that the basic transformations in Definition 2 are restricted versions of them.
\begin{theorem}
If the stopping model, $\mathcal{N} = \{N_{\eta}: h_{N_{\eta}}, ~ \eta \in [\eta_0,\infty) \}$, is
``\emph{uniparametric and closed under pgf composition}" as defined above, then:
$$  h_{N_{\eta_2}}^{-1}\circ h_{N_{\eta_1}} = h_{N_{\eta_1}}\circ h_{N_{\eta_2}}^{-1}  $$
for all $\eta_1, \eta_2 \in [\eta_0,\infty)$. Furthermore, $H_{\mathcal{N}}(t;\eta_1,\eta_2) =
h_{N_{\eta_{1}}}\circ h_{N_{\eta_{2}}}^{-1}(t)$
can be parametrized in an identifiable way through
$\eta = - \log(H_{\mathcal{N}}^{'}(0;\eta_1,\eta_2)) = \eta_1 - \eta_2$, and if one denotes
$H_{\mathcal{N},\eta}(t) = H_{\mathcal{N}}(t;\eta_1,\eta_2)$ with $\eta \in \mathbb{R}$, then
\begin{enumerate}\vspace{-.3cm}
\item{} $H_{\mathcal{N},\eta}\circ H_{\mathcal{N},\eta^{'}} = H_{\mathcal{N},\eta + \eta^{'}}$ for all $\eta,\eta^{'} \in \mathbb{R}$,
\item{} when $\eta \ge \eta_0$, then $H_{\mathcal{N},\eta} = h_{N_{\eta}}$,
\item{} when $\eta \ge \eta_0$, then $H_{\mathcal{N},-\eta} = h_{N_{\eta}}^{-1}$,
\item{} $H_{\mathcal{N},\eta=0}(t) = t$.
\end{enumerate}\vspace{-.3cm}
Likewise,
$ \overline{h}_{N_{\eta_1}}\circ \overline{h}_{N_{\eta_2}}^{-1} =
\overline{h}_{N_{\eta_2}}^{-1}\circ \overline{h}_{N_{\eta_1}},$
and the properties listed above also apply for
$\overline{H}_{\mathcal{N},\eta}(t) = \overline{H}_{\mathcal{N}}(t;\eta_1,\eta_2) = \overline{h}_{N_{\eta_1}}\circ \overline{h}_{N_{\eta_2}}^{-1}(t)
= 1 - H_{\mathcal{N,\eta}}(1-t).$
\end{theorem}
{\sl{Proof:}} The commutativity for $\eta_1, \eta_2 \in [\eta_0,\infty)$ follows from:
$$h_{N_{\eta_2}}^{-1}\circ h_{N_{\eta_1}} =
h_{N_{\eta_2}}^{-1}\circ h_{N_{\eta_1}}\circ h_{N_{\eta_2}}\circ h_{N_{\eta_2}}^{-1} =
h_{N_{\eta_2}}^{-1}\circ h_{N_{\eta_2}}\circ h_{N_{\eta_1}}\circ h_{N_{\eta_2}}^{-1} =
h_{N_{\eta_1}}\circ h_{N_{\eta_2}}^{-1}.$$
$H_{\mathcal{N}}(t;\eta_1,\eta_2)$ can be parametrized through
$\eta = - \log(H_{\mathcal{N}}^{'}(0;\eta_1,\eta_2)) = \eta_1 - \eta_2$ because
$$H^{\prime}_{\mathcal{N}}(0;\eta_1,\eta_2) =
\left(h_{N_{\eta_2}}^{-1}\right)^{\prime}(0)\cdot h_{N_{\eta_1}}^{\prime}(0)=
\frac{1}{h_{N_{\eta_2}}^{\prime}(0)}h_{N_{\eta_1}}^{\prime}(0)=
e^{-(\eta_1-\eta_2)},$$
and if $\eta=\eta_1 - \eta_2 = \eta_1^{\prime} - \eta_2^{\prime}$ with $\eta_1,\eta_2,\eta_1^{\prime},\eta_2^{\prime} \ge \eta_{0}$, then:
$$h_{N_{\eta_2}}^{-1}\circ h_{N_{\eta_1}} =
h_{N_{\eta_2^{\prime}}}^{-1}\circ h_{N_{\eta_2^{\prime}}}\circ h_{N_{\eta_2}}^{-1}\circ h_{N_{\eta_1}} =
h_{N_{\eta_2^{\prime}}}^{-1}\circ h_{N_{\eta_2^{\prime}}}\circ h_{N_{\eta_1}}\circ h_{N_{\eta_2}}^{-1} =
h_{N_{\eta_2^{\prime}}}^{-1}\circ h_{N_{\eta_1}+N_{\eta_2^{\prime}}} \circ h_{N_{\eta_2}}^{-1} = $$
$$h_{N_{\eta_2^{\prime}}}^{-1}\circ h_{N_{\eta_1^{\prime}}+N_{\eta_2}} \circ h_{N_{\eta_2}}^{-1} =
h_{N_{\eta_2^{\prime}}}^{-1}\circ h_{N_{\eta_1^{\prime}}}\circ h_{N_{\eta_2}}\circ h_{N_{\eta_2}}^{-1} =
h_{N_{\eta_2^{\prime}}}^{-1}\circ h_{N_{\eta_1^{\prime}}},$$
and because if $h_{N_{\eta_2}}^{-1}\circ h_{N_{\eta_1}} = h_{N_{\eta_2^{\prime}}}^{-1}\circ h_{N_{\eta_1^{\prime}}}$
with $\eta_1,\eta_2,\eta_1^{\prime},\eta_2^{\prime} \ge \eta_{0}$, then:
$$ \left(h_{N_{\eta_2}}^{-1}\circ h_{N_{\eta_1}}\right)_{|t=0}^{\prime}=
\left(h_{N_{\eta_2^{\prime}}}^{-1}\circ h_{N_{\eta_1^{\prime}}}\right)_{|t=0}^{\prime}, $$
and so $e^{-(\eta_1-\eta_2)}= e^{-(\eta_1^{\prime}-\eta_2^{\prime})}$ and $\eta_1-\eta_2 = \eta_1^{\prime} - \eta_2^{\prime}.$
To prove the additivity of $H_{\mathcal{N},\eta}\circ H_{\mathcal{N},\eta^{'}}$, let $\beta \ge \eta_{0}+ \max(|\eta|,|\eta^{\prime}|)$ and note that:
$$ H_{\mathcal{N},\eta}\circ H_{\mathcal{N},\eta^{\prime}} =
(h_{N_{\beta}}^{-1}\circ h_{N_{\beta+\eta}})\circ (h_{N_{\beta}}^{-1}\circ h_{N_{\beta+\eta^{\prime}}}) = $$
$$ h_{N_{\beta}}^{-1}\circ h_{N_{\beta}}^{-1}\circ h_{N_{\beta+\eta}}\circ h_{N_{\beta+\eta^{\prime}}} =
h_{N_{2 \beta}}^{-1}\circ h_{N_{2 \beta +\eta+\eta^{\prime}}} = H_{\mathcal{N},\eta+\eta^{\prime}}. $$
Furthermore, letting $\beta \ge \eta_0$ one has that for any $\eta \ge \eta_0$:
$$H_{\mathcal{N},\eta} = h_{N_{\beta}}^{-1}\circ h_{N_{\beta+\eta}} =
h_{N_{\beta}}^{-1}\circ h_{N_{\beta}}\circ h_{N_{\eta}}=h_{N_{\eta}}, $$
$$H_{\mathcal{N},-\eta} = h_{N_{\beta+\eta}}^{-1}\circ h_{N_{\beta}} = h_{N_{\eta}}^{-1}\circ h_{N_{\beta}}^{-1}\circ h_{N_{\beta}} = h_{N_{\eta}}^{-1},$$
and that, $H_{\mathcal{N},\eta=0}(t) = h_{N_{\beta}}^{-1}\circ h_{N_{\beta}}(t) = t$.
\hfill{$\square$}

By using $H_{\mathcal{N},\eta}$ or $\overline{H}_{\mathcal{N},\eta}$ with $\eta \in \mathbb{R}$ in a model extension
of Definition 3, one extends the parameter space through
values of $\eta$ in the whole real line and not just in $\mathcal{H} = [\eta_0,\infty)$.

Many stopping models satisfy the consequences of Theorem 1 without being closed under pgf composition.
Next, an extra necessary condition for being a stopping model closed under pgf composition
is obtained by imposing that the $t^2$ coefficients of the series expansion of $h_{N_{\eta_{2}}}^{-1}\circ h_{N_{\eta_{1}}}$
and of $h_{N_{\eta_{1}}}\circ h_{N_{\eta_{2}}}^{-1}$ have to be equal for any $\eta_1, \eta_2 \in [\eta_0,\infty)$.
Imposing that higher order term coefficients of these expansions are equal leads to other necessary conditions.
\begin{corollary}
If a stopping model $\mathcal{N} = \{N_{\eta}: h_{N_{\eta}}, ~ \eta \in [\eta_0,\infty) \}$
is closed under pgf composition,
$$ \frac{\Pr(N_{\eta}=2)}{\Pr(N_{\eta}=1)(1-\Pr(N_{\eta}=1))} = C, ~~ \mbox{for all} ~~ \eta \in [\eta_0,\infty). $$
\end{corollary}

\subsection{Examples of stopping models closed under pgf composition}

\textbf{Example 6.1:} The \emph{potential conjugate} model,
$$ \mathcal{N} = \{ N_p: h_{N_p} = 1 - (1 - t)^p, ~~ p \in (0,1]  \}, $$
is closed under pgf composition with $E[N_p]=\infty$ and $\Pr(N_p=1)=p$, and therefore with $\eta=-\log{p} \in [0,\infty)$.
It includes $N_I$ but it is not auto-reversible as described in Section 5.3.

\textbf{Example 6.2:} The zero-truncated geometric model,
$$ \mathcal{N} = \{ N_p: h_{N_p} = \frac{pt}{1-(1-p)t}, ~~ p \in (0,1] \},  $$
is closed under pgf composition, with
$\Pr(N_p=1) = p$ and $\eta = - \log{p} \in [0,\infty)$.
It includes $N_I$ and, as indicated in Example 5.4, it is auto-reversible.

The next result provides a way of generating a new model closed under pgf composition, starting from an initial model closed
under pgf composition and two random variables that do not belong to the initial model but whose pgf composition does.
\begin{proposition}
Let the stopping model $\mathcal{N} = \{N_{\eta}: h_{N_{\eta}}, ~ \eta \in [\eta_0,\infty) \}$
be closed under pgf composition,
and let $N_1$ and $N_2$ be two random variables that do not belong to $\mathcal{N}$ but
such that $h_{N_1}\circ h_{N_2} = h_{N_{\alpha}}$ with $N_{\alpha} \in \mathcal{N}$.
Then, for any given $\alpha > 0$ the statistical model
$$\mathcal{N}_{\alpha} = \{ \tilde{N}_{\eta} : h_{\tilde{N}_{\eta}} = h_{N_2}\circ h_{N_{\eta-\alpha}}\circ h_{N_1}, ~~ \eta \in [\alpha+\eta_0,\infty) \}      $$
is also closed under pgf composition.
\end{proposition}
{\sl{Proof:}} If $\tilde{N}_{\eta_1}$ and $\tilde{N}_{\eta_2}$ are random variables that belong to $\mathcal{N}_{\alpha}$, then
$$   h_{\tilde{N}_{\eta_1}} \circ h_{\tilde{N}_{\eta_2}} =  h_{N_2} \circ h_{N_{\eta_1-\alpha}} \circ h_{N_1} \circ h_{N_2} \circ h_{N_{\eta_2-\alpha}} \circ h_{N_1} =  $$
$$   h_{N_2} \circ h_{N_{\eta_1-\alpha}} \circ h_{N_{\alpha}} \circ h_{N_{\eta_2-\alpha}} \circ h_{N_1}  =  h_{N_2} \circ h_{N_{\eta_1+\eta_2-\alpha}} \circ h_{N_1},          $$
which is the pgf of a random variable  $\tilde{N}_{\eta_1+\eta_2}$ that also belongs to $\mathcal{N}_{\alpha}$.
\hfill{$\square$}

Using Proposition 3 twice in a row does not generate any new family of models.
Next, this result is used to generate three families of stopping models closed under
pgf composition starting from the geometric model.

\textbf{Example 6.3:}
If $N_1$ is zero-truncated Poisson($\alpha$)
and $N_2$ is Logarithmic($\alpha$),
as in Example 5.2, then $h_{N_1}\circ h_{N_2}$ is the pgf of a Geometric($p=e^{-\alpha}$) and by
Proposition 3 one has that for any given value of $\alpha > 0$ the statistical model:
$$ \mathcal{N}_{\alpha} = \{ N_{\eta} : h_{N_{\eta}} = \frac{1}{\alpha} \ln\left(1 + \frac{(e^{\alpha t} -1)(e^{\alpha} - 1)}
{(e^{\eta} -1)(e^{\alpha} - e^{\alpha t}) + e^{\alpha} -1}\right), ~~ \eta \in [\alpha,\infty) \}, $$
is closed under pgf composition. In the limit, when $\alpha$ tends to $0$ this model becomes the
geometric model, and when $\alpha$ tends to $\infty$ it becomes $N_I$.
The model $\mathcal{N}_{\alpha}$ is extreme auto-reversible for every $\alpha$, but it only includes $N_I$ in the limiting cases mentioned.

\textbf{Example 6.4:}
Let $N_1$ be zero-truncated Negative-Binomial($\alpha \beta, \beta$), with:
$$ h_{N_1} = \frac{(1-(1-e^{-\frac{\alpha}{\beta}})t)^{-\beta}-1}{e^{\alpha}-1}, $$
and let $N_2$ be extended truncated Negative-Binomial($\alpha,-1/\beta$) in Engen (1974), with:
$$ h_{N_2} = \frac{(1-(1-e^{-\alpha})t)^{\frac{1}{\beta}}-1}{e^{-\frac{\alpha}{\beta}}-1}, $$
where $\alpha \ge 0$ and $\beta \ge 1$. Then $h_{N_1}\circ h_{N_2}$ is the pgf of a Geometric($p=e^{-\alpha}$),
and by Proposition 3 one has that given any $\alpha \ge 0$ and $\beta \ge 1$ the statistical model:
$$ \mathcal{N}_{\alpha,\beta} = \{ N_{\eta} : h_{N_{\eta}} = \frac{1-\left(\frac{(-{{\rm e}^{\alpha+\eta}}+1)(1-t+t{{\rm e}^{-\frac{\alpha}{\beta}}})^{\beta}+{{\rm e}^{\eta}}-1}{({{\rm e}^{\alpha}}-{{\rm e}^{\alpha+\eta}})(1-t+t{{\rm e}^{-\frac{\alpha}{\beta}}})^{\beta}+{{\rm e}^{\eta}}-{{\rm e}^{\alpha}}}\right)^{\frac{1}{\beta}}}{1-{\rm e}^{-\frac{\alpha}{\beta}}}, ~~ \eta \in [\alpha,\infty) \},  $$
is closed under pgf composition.
When $\beta$ tends to $\infty$ one obtains the models in Example 6.3, and when $\alpha$ tends to $0$ or
$\beta$ converges to $1$ one obtains the geometric model in Example 6.2.
Other than in these limiting cases, $\mathcal{N}_{\alpha,\beta}$ is neither extreme auto-reversible, nor includes $N_I$.

\textbf{Example 6.5:}
Let $N_1$ be zero-truncated Binomial($n, p = 1 - e^{-\alpha/n}$), with:
$$ h_{N_1} = \frac{(1+(e^{\frac{\alpha}{n}}-1)t)^{n}-1}{e^{\alpha}-1}, $$
and let $N_2$ be zero-truncated Negative-Binomial($\alpha,1/n$), with:
$$ h_{N_2} = \frac{(1-(1-e^{-\alpha})t)^{-\frac{1}{n}}-1}{e^{\frac{\alpha}{n}}-1}, $$
where $\alpha \ge 0$ and $n \in \mathbb{N}^{+}$. Then, $h_{N_1}\circ h_{N_2}$ is the pgf of a Geometric($p=e^{-\alpha}$), and
by Proposition 3 one has that for any given $\alpha \ge 0$ and $n \in \mathbb{N}^{+}$ the statistical model
$$  \mathcal{N}_{\alpha,n} = \{ N_{\eta} : h_{N_{\eta}} =
\frac{\left(\frac{({{\rm e}^{\eta}}-1)(-t+1+t{{\rm e}^{{\frac{\alpha}{n}}}})^{n}-{{\rm e}^{\alpha+\eta}}+1}{(-{{\rm e}^{\alpha}}+{{\rm e}^{\eta}})(-t+1+t{{\rm e}^{{\frac{\alpha}{n}}}})^{n}+{{\rm e}^{\alpha}}-{{\rm e}^{\alpha+\eta}}}\right)^{-\frac{1}{n}}-1}{{{\rm e}^{{\frac{\alpha}{n}}}}-1}, ~~ \eta \in [\alpha,\infty) \},  $$
is closed under pgf composition. In the limit,
when $n$ converges to $\infty$ one obtains the models in Example 6.3, and when $\alpha$ converges to $0$, or
when $n$ is $1$, one obtains the geometric model in Example 6.2.
Other than in these limiting cases, $\mathcal{N}_{\alpha,n}$ is neither extreme auto-reversible, nor includes $N_I$,
but it is extreme reversible with the $\mathcal{N}_{\alpha,\beta=n}$ in Example 6.4.

The next result
provides a way of generating a family of statistical models
closed under pgf composition starting from any model that is like that.
\begin{proposition}
If the stopping model $\mathcal{N} = \{N_{\eta}: h_{N_{\eta}}(t), ~ \eta \in [\eta_0,\infty) \}$
is closed under pgf composition then, for every given $k \in \mathbb{N}^{+},$ the statistical model
$$\mathcal{N}_{k} = \{ \tilde{N}_{\eta} : h_{\tilde{N}_{\eta}}(t) = \left(h_{N_{\eta}}(t^k)\right)^{1/k}, ~~ \eta \in [\eta_0,\infty) \}      $$
is also closed under pgf composition.
\end{proposition}
Using this result with Examples 6.1 and 6.2 one has that for every $k \in \mathbb{N}^{+}$
the models
$$ \mathcal{N}_k = \{ N_p: h_{N_p} = \left(1 - \left(1 - t^k\right)^p\right)^{1/k}, ~~ p \in (0,1]  \}, $$
and
$$ \mathcal{N}_k = \{ N_p: h_{N_p} = \left(\frac{pt^k}{1-(1-p)t^k} \right)^{1/k}, ~~ p \in (0,1] \},  $$
are closed under pgf composition with $\eta = - (1/k) \log{p}$ and support $n = 1,k+1,2k+1,\ldots$.

Finally we present a family of statistical models closed under pgf composition that embed Examples 6.1 and 6.2
as limiting cases and all include $N_I$.

\textbf{Example 6.6:}
Given any value $\alpha \in (0,1)$, the statistical model
$$ \mathcal{N}_{\alpha} = \{ N_p: h_{N_p} = 1 - \frac{1-t}{(p+(1-p)(1-t)^\alpha)^{1/\alpha}}, ~~ p \in (0,1]  \}, $$
is closed under pgf composition with $E[N_p]=p^{-1/\alpha}$, with $Var[N_p]=\infty$ and with $\eta=-\log{p} \in [0,\infty)$,
and it contains $N_I$. In the limit, when $\alpha$ tends to $0$ one obtains the model in Example 6.1, and when $\alpha$ tends to $1$
one obtains the model in Example 6.2.

\section{Randomly stopped extreme based model transformations}

\subsection{Model transformations in Definitions 2 and 3}

Given the properties of $h_N$ and of $\overline{h}_N^{-1}$ it follows that
$F_{max_N(X)}(y) \le F_X(y)$ and $F_{min_N^{-1}(X)}(y) \le F_X(y)$ for all $y$ in their domain, and therefore that
$max_N(X)$ and $min_N^{-1}(X)$ are random variables always larger than $X$ in the usual stochastic order.
Furthermore, given the properties of $\overline{h}_N$ and of $h_N^{-1}$ it follows that
$F_{min_N(X)}(y) \ge F_X(y)$ and $F_{max_N^{-1}(X)}(y) \ge F_X(y)$, and therefore that
$min_N(X)$ and $max_N^{-1}(X)$ are always smaller than $X$ in that stochastic order.

Hence, two of the basic transformations of Definition 2 transform any model $\mathcal{X} = \{X_{\theta}, ~ \theta \in \Theta \}$
into a model $\mathcal{Y} = \{Y_{\theta,\delta}, ~ \theta \in \Theta, \delta \in \mathcal{D} \}$
with random variables $Y_{\theta,\delta}$ stochastically larger than $X_{\theta}$, while the other two
transform $\mathcal{X}$ into a model with $Y_{\theta,\delta}$ stochastically smaller than $X_{\theta}$.

The four combined transformations of Definition 3 transform
$\mathcal{X}$ into a model $\mathcal{Y}$ with random variables $Y_{\theta,\delta_1,\delta_2}$ that can be
stochastically larger and smaller than $X_{\theta}$.

By construction,
the dimension of the parameter space of models obtained through
transformations in Definition 3 is never smaller than the dimension of the parameter space of models obtained through
transformations in Definition 2, which in turn is never smaller than the dimension of the parameter space of the initial model.
We next investigate when is the initial model always included in the transformed model,
and when does repeated use of these extensions leave the extended model unchanged.

\subsection{When do transformations work as extensions?}

A sufficient condition for basic transformations in Definition 2 to work as extensions for any model, $\mathcal{X}$,
is that the identity belongs to the stopping model.
\begin{proposition}
If $N_I \in \mathcal{N}$, with $\Pr(N_I=1)=1$, then the four basic model transformations in Definition 2
work as a model extension of $\mathcal{X}$, for any $\mathcal{X}$.
\end{proposition}
{\sl{Proof:}} If $N_I$, with $h_{N_I}(t) = t$, belongs to $\mathcal{N}$,
then $X \in \mathcal{X}$ implies that $X \in max_{\mathcal{N}}(\mathcal{X})$, and so
$\mathcal{X} \subset max_{\mathcal{N}}(\mathcal{X})$. The same argument applies to the other
three basic transformations.
\hfill{$\square$}

If one starts with a single random variable, $\mathcal{X} = \{X\}$,
then $N_I \in \mathcal{N}$ is necessary and sufficient for $X$ to be included in
$max_{\mathcal{N}}(X)$ and in $min_{\mathcal{N}}(X)$. In general though, one can find instances of specific models,
$\mathcal{X}$, included in $max_{\mathcal{N}}(\mathcal{X})$ or in $min_{\mathcal{N}}(\mathcal{X})$ without $N_I$
belonging to $\mathcal{N}$.

On the other hand, the four combined mechanisms of Definition 3 always work as model extensions,
irrespective of whether $N_I$ is in $\mathcal{N}$ or not.
\begin{proposition}
The four model transformation in Definition 3
work as a model extension of $\mathcal{X}$, for any $\mathcal{X}$. That is so, even if one of the two new
parameters, $\delta_1$ or $\delta_2$, is fixed.
\end{proposition}
{\sl{Proof:}}
$F_{X_{\theta}} \in \mathcal{X}$ implies that
$F_{Y_{\theta,\delta_1,\delta_2}} = h_{N_{\delta_2}}\circ h_{N_{\delta_1}}^{-1}(F_{X_{\theta}}) \in max_{\mathcal{N}}(max^{-1}_{\mathcal{N}}(\mathcal{X}))$
for all $\delta_1,\delta_2 \in \mathcal{D}$,
and in particular $F_{X_{\theta}} = h_{N_{\delta_1}}\circ h_{N_{\delta_1}}^{-1}(F_{X_{\theta}}) \in max_{\mathcal{N}}(max^{-1}_{\mathcal{N}}(\mathcal{X}))$,
which means that $\mathcal{X} \subset max_{\mathcal{N}}(max^{-1}_{\mathcal{N}}(\mathcal{X}))$. The same argument applies to the
other three transformations, and when any of the two new parameters is fixed.
\hfill{$\square$}

Different from the transformations in Definition 3, using
$max_{\mathcal{N}}(\cdot) \cup min_{\mathcal{N}}(\cdot)$ with a stopping model $\mathcal{N}$ that does not include $N_I$
does not always work as a model extension.

\subsection{When are the extensions statistically stable?}

Under general uniparametric stopping models, the basic transformations of Definition 2
usually add one dimension to the parameter space,
while the combined transformations of Definition 3 usually add two dimensions to it.

Instead, when the stopping model is uniparametric and closed under pgf composition
both basic as well as combined transformations add at most a single dimension,
and the basic transformations of Definition 2 become restricted versions of the combined transformations
of Definition 3 with the extra parameter, $\eta$, of the basic transformations
taking values on a semi-line and the extra parameter, $\eta$, of the combined transformation taking
values on the whole real line.

Furthermore, under general stopping models repeated use of these extensions usually keep extending the models.
Instead, when the stopping model is closed under pgf composition and the transformation works as an extension,
then it is always a statistically stable extension and hence repeated use of that extension leaves the extended
model unchanged.
\begin{proposition}
If the stopping model $\mathcal{N} = \{N_{\eta}: h_{N_{\eta}}, ~ \eta \in [\eta_0,\infty) \}$
is ``\emph{uniparametric and closed under pgf composition}" then,
\begin{enumerate}\vspace{-.3cm}
\item{}
if $\eta_0 =0$, then
$max_{\mathcal{N}}(\cdot)$,
$min_{\mathcal{N}}(\cdot)$,
$max_{\mathcal{N}}^{-1}(\cdot)$,
and $min_{\mathcal{N}}^{-1}(\cdot)$
are statistical model extensions that are statistically stable, and
\item{}
if $\eta_0 > 0$, then
$max_{\mathcal{N}}(\mathcal{X})$ is contracted by $max_{\mathcal{N}}(\cdot)$,
$min_{\mathcal{N}}(\mathcal{X})$ is contracted by $min_{\mathcal{N}}(\cdot)$,
$max_{\mathcal{N}}^{-1}(\mathcal{X})$ is contracted by $max_{\mathcal{N}}^{-1}(\cdot)$, and
$min_{\mathcal{N}}^{-1}(\mathcal{X})$ is contracted by $min_{\mathcal{N}}^{-1}(\cdot)$, for all $\mathcal{X}$.
\end{enumerate}\vspace{-.3cm}
\end{proposition}
{\sl{Proof:}}
By Theorem 2 one has that for any $\mathcal{X}$,
$$ \mathcal{Y} = max_{\mathcal{N}}(max_{\mathcal{N}}(\mathcal{X})) = \{Y_{\theta,\eta_1,\eta_2}: F_{Y_{\theta,\eta_1,\eta_2}} =
h_{N_{\eta_2}}\circ h_{N_{\eta_1}}(F_{X_{\theta}}), ~ \theta \in \Theta, \eta_1,\eta_2 \in [\eta_0,\infty) \} = $$
$$\{Y_{\theta,\eta}: F_{Y_{\theta,\eta}} = h_{N_{\eta=\eta_1+\eta_2}}(F_{X_{\theta}}), ~ \theta \in \Theta, \eta \in [2 \eta_0,\infty) \} \subset  $$
$$ \{Y_{\theta,\eta}: F_{Y_{\theta,\eta}} = h_{N_{\eta}}(F_{X_{\theta}}), ~ \theta \in \Theta, \eta \in [\eta_0,\infty) \} = max_{\mathcal{N}}(\mathcal{X}),  $$
and so if $\eta_0 > 0$, then $max_{\mathcal{N}}(\cdot)$ contracts $max_{\mathcal{N}}(\mathcal{X})$.
When $\eta_0 = 0$,
$$  max_{\mathcal{N}}(max_{\mathcal{N}}(\mathcal{X}))  =
\{Y_{\theta,\eta}: F_{Y_{\theta,\eta}} = h_{N_{\eta}}(F_{X_{\theta}}), ~ \theta \in \Theta, \eta \in [0,\infty) \} =  max_{\mathcal{N}}(\mathcal{X}), $$
which means that
$max_{\mathcal{N}}(\cdot)$ is a statistically stable extension.
The same argument applies to the other three transformations in Definition 2.
\hfill{$\square$}

The next result establishes that under uniparametric stopping models closed under pgf composition,
there are only two distinct combined extensions and they are statistically stable.
\begin{proposition}
If the stopping model $\mathcal{N}$ is ``\emph{uniparametric and closed under pgf composition}", then Definition 3 yields only two distinct
model extensions which are:
\begin{enumerate}\vspace{-.3cm}
\item{} $\mathcal{Y}_1 = max^{-1}_{\mathcal{N}}(max_{\mathcal{N}}(\cdot)) = max_{\mathcal{N}}(max_{\mathcal{N}}^{-1}(\cdot))$, and
\item{} $\mathcal{Y}_2 = min^{-1}_{\mathcal{N}}(min_{\mathcal{N}}(\cdot)) = min_{\mathcal{N}}(min_{\mathcal{N}}^{-1}(\cdot))$,
\end{enumerate}\vspace{-.3cm}
and these two extensions are both statistically stable. Furthermore, in that case it holds that:
\begin{enumerate}\vspace{-.3cm}
\item{} the model $\mathcal{Y}_1 = max^{-1}_{\mathcal{N}}(max_{\mathcal{N}}(\mathcal{X}))$ is invariant under $max_{\mathcal{N}}(\cdot)$
and $max_{\mathcal{N}}^{-1}(\cdot)$,
\item{} the model $\mathcal{Y}_2 = min^{-1}_{\mathcal{N}}(min_{\mathcal{N}}(\mathcal{X}))$ is invariant under $min_{\mathcal{N}}(\cdot)$
and $min_{\mathcal{N}}^{-1}(\cdot)$, and
\item{} the transformations in Definition 2 are restricted versions of one of these two extensions.
\end{enumerate}\vspace{-.3cm}
\end{proposition}
{\sl{Proof:}}
By Theorem 3, one has that for any $\mathcal{X}$,
$$  \mathcal{Y}_1 = max_{\mathcal{N}}^{-1}(max_{\mathcal{N}}(\mathcal{X})) = max_{\mathcal{N}}(max_{\mathcal{N}}^{-1}(\mathcal{X}))  =
\{Y_{\theta,\eta}: F_{Y_{\theta,\eta}} = H_{\mathcal{N},\eta}(F_{X_{\theta}}), ~ \theta \in \Theta, \eta \in {\mathcal{R}} \},$$
and the stability follows from that same theorem, because:
$$ max_{\mathcal{N}}(\mathcal{Y}_1) =
\{Y_{\theta,\eta+\eta^{'}}: F_{Y_{\theta,\eta+\eta^{'}}} = H_{\mathcal{N},\eta+\eta^{'}}(F_{X_{\theta}}), ~ \theta \in \Theta, \eta+\eta^{'} \in {\mathcal{R}} \} =
\mathcal{Y}_1, $$
$$ max^{-1}_{\mathcal{N}}(\mathcal{Y}_1) =
\{Y_{\theta,\eta-\eta^{'}}: F_{Y_{\theta,\eta-\eta^{'}}} = H_{\mathcal{N},\eta-\eta^{'}}(F_{X_{\theta}}), ~ \theta \in \Theta, \eta-\eta^{'} \in {\mathcal{R}} \} =
\mathcal{Y}_1. $$
By Theorem 3 one also has that:
$$  \mathcal{Y}_2 = min_{\mathcal{N}}^{-1}(min_{\mathcal{N}}(\mathcal{X})) =  min_{\mathcal{N}}(min_{\mathcal{N}}^{-1}(\mathcal{X})) =
\{Y_{\theta,\eta}: F_{Y_{\theta,\eta}} = \overline{H}_{\mathcal{N,\eta}}(F_{X_{\theta}}), ~ \theta \in \Theta, \eta \in {\mathcal{R}} \},$$
where $\overline{H}_{\mathcal{N,\eta}}(t) = 1 - H_{\mathcal{N,\eta}}(1-t)$, and stability follows likewise.
\hfill{$\square$}

\begin{corollary}
If $\mathcal{N}$ is ``\emph{uniparametric and closed under pgf composition}", then:
\begin{enumerate}\vspace{-.3cm}
\item{} $max_{\mathcal{N}}(\mathcal{X}) \cup max_{\mathcal{N}}^{-1}(\mathcal{X}) \subset
\mathcal{Y}_1 = max^{-1}_{\mathcal{N}}(max_{\mathcal{N}}(\mathcal{X})) =
max_{\mathcal{N}}(max^{-1}_{\mathcal{N}}(\mathcal{X}))$,
\item{} $min_{\mathcal{N}}(\mathcal{X}) \cup min_{\mathcal{N}}^{-1}(\mathcal{X}) \subset
\mathcal{Y}_2 = min^{-1}_{\mathcal{N}}(min_{\mathcal{N}}(\mathcal{X})) =
min_{\mathcal{N}}(min^{-1}_{\mathcal{N}}(\mathcal{X}))$,
\end{enumerate}\vspace{-.3cm}
and if $N_I \in \mathcal{N}$,
then the models on the left and the models on the right are equal.
\end{corollary}

\subsection{What happens with stopping models both closed and extreme reversible?}

When two stopping models are closed under pgf composition and extreme reversible,
Proposition 8 and the definition of extreme reversibility lead to the next result.
\begin{proposition}
If the stopping models $\mathcal{N}$ and $\mathcal{N}^{*}$ are uniparametric, closed under pgf composition, and extreme reversible, then
the two distinct statistically stable model extensions in Definition 3 obtained with $\mathcal{N}$
and the ones obtained with $\mathcal{N}^{*}$ are the same extensions.
\end{proposition}
According to Proposition 8, when a stopping model is closed under pgf composition the four
extensions in Definition 3 collapse down into two distinct ones.
The next result, stating that when a stopping model is both closed and extreme auto-reversible
then these two extensions become a single one, is a straight consequence of the definition of extreme-auto-reversibility.
\begin{proposition}
If the stopping model $\mathcal{N}$ is uniparametric, closed under pgf composition, and extreme auto-reversible, then
the four statistically stable model extensions in Definition 3 coincide,
$$ max^{-1}_{\mathcal{N}}(max_{\mathcal{N}}(\cdot)) = max_{\mathcal{N}}(max_{\mathcal{N}}^{-1}(\cdot)) =
min^{-1}_{\mathcal{N}}(min_{\mathcal{N}}(\cdot)) = min_{\mathcal{N}}(min_{\mathcal{N}}^{-1}(\cdot)), $$
and they coincide with $min_{\mathcal{N}}(max_{\mathcal{N}}(\cdot))$ and with
$max_{\mathcal{N}}(min_{\mathcal{N}}(\cdot))$.
If on top of that, $N_I \in \mathcal{N}$ (i.e. $\eta_0=0$),
this statistically stable model extension also coincides with $max_{\mathcal{N}}(\cdot) \cup min_{\mathcal{N}}(\cdot)$.
\end{proposition}
The geometric stopping model satisfies all the conditions of Proposition 10. As a consequence,
the Marshall-Olkin extension of $\mathcal{X}$, originally defined
to be $max_{\mathcal{N}}(\mathcal{X}) \cup min_{\mathcal{N}}(\mathcal{X})$
when $\mathcal{N}$ is geometric, coincides with the extension of $\mathcal{X}$ obtained through Definition 3 with geometric stopping.

As a consequence, we consider the model extensions in Definition 3 to be the natural
way to generalize the Marshall-Olkin extension with stopping models other than geometric.
Different from what happens if one generalized Marshall-Olkin through $max_{\mathcal{N}}(\cdot) \cup min_{\mathcal{N}}(\cdot)$,
by generalizing them through the transformations in Definition 3 one guarantees that these transformations will
work as model extensions under any stopping model, $\mathcal{N}$.

\section{Examples of statistically stable extensions}

When one uses the model extensions of Definition 3 with stopping models that are neither closed under pgf
composition nor extreme auto-reversible, one obtains four different extensions that are not statistically stable, and
the four basic transformations of Definition 2 are not restricted versions of them. As an example, Appendix 1 presents the four
basic and the four combined extensions obtained when $\mathcal{N}$ are the zero-truncated Poisson or the logarithmic models.

Here we present the model extensions in Definition 3 obtained when the
stopping models are the ones presented in Section 6.2. Given that these stopping models are all closed under pgf composition,
all the extensions obtained here are statistically stable in the sense that applying them twice on any given
model leads to the same model as applying them once.

Furthermore, because of Proposition 8 another consequence of
all these stopping models being closed under pgf composition
is that for them Definition 3 yields at most two
distinct extensions, and that the transformations in Definition 2
are restrictions of these two extensions and do not need to be considered apart.

In three of the examples, the stopping models are not auto-reversible, and for them the model extension in Definition 3 based on
maxima extends $\mathcal{X} = \{X_{\theta}: F_{X_{\theta}}, ~ \theta \in \Theta \}$
into:
$$  \mathcal{Y}_1 = max_{\mathcal{N}}(max^{-1}_{\mathcal{N}}(\mathcal{X})) =
 \{Y_{\theta,\eta}: F_{Y_{\theta,\eta}} = H_{\mathcal{N},\eta}(F_{X_{\theta}}), ~~ \theta \in \Theta, \eta \in (-\infty,\infty) \},  $$
and the extension in Definition 3 based on minima extends $\mathcal{X}$ into:
$$  \mathcal{Y}_2 =  min_{\mathcal{N}}(min^{-1}_{\mathcal{N}}(\mathcal{X})) =
 \{Y_{\theta,\eta}: F_{Y_{\theta,\eta}} = \overline{H}_{\mathcal{N},\eta}(F_{X_{\theta}}), ~~ \theta \in \Theta, \eta \in (-\infty,\infty) \}, $$
with $H_{\mathcal{N},\eta}(\cdot)$ and $\overline{H}_{\mathcal{N},\eta}(\cdot)$ as in Theorem 3.

In the second and third examples the stopping models are auto-reversible, and hence for them
these two extension mechanisms, $\mathcal{Y}_1$ and $\mathcal{Y}_2$, coincide because of Proposition 10.
The fourth and fifth families of stopping models considered can be reversible, and when they are reversible they
lead to the same pair of model extensions because of Proposition 9.

\textbf{Example 8.1:} Let the stopping model be the one in Example 6.1,
$$ \mathcal{N} = \{N_{\eta}: h_{N_{\eta}} = 1 - (1 - t)^{e^{-\eta}}, ~ \eta \in [0,\infty) \},  $$
which is not extreme auto-reversible but it includes $N_I$.

The extension of $\mathcal{X}$ obtained through maxima and precursors of maxima is:
$$  \mathcal{Y}_1 = max_{\mathcal{N}}(max^{-1}_{\mathcal{N}}(\mathcal{X})) = \{Y_{\theta,\eta}: F_{Y_{\theta,\eta}} = 1-(1-F_{X_{\theta}})^{e^{-\eta}}, ~ \theta \in \Theta, \eta \in (-\infty,\infty) \},  $$
which is a special case of the extension in Cordeiro and Castro (2009).
When one restricts $\eta \ge 0$, here one obtains $\mathcal{Y}_1^{\prime} = max_{\mathcal{N}}(\mathcal{X})$,
and when one restricts $\eta \le 0$ one obtains
$\mathcal{Y}_1^{\prime \prime} = max^{-1}_{\mathcal{N}}(\mathcal{X})$, and therefore in this case
$\mathcal{Y}_1 = max_{\mathcal{N}}(\mathcal{X}) \cup max^{-1}_{\mathcal{N}}(\mathcal{X})$.

When $\mathcal{X}$ is for example an exponential random variable, $\mathcal{Y}_1$ becomes the exponential model.
Because of the stability of this extension,
using it again, now on the exponential model, will leave that model unchanged which means that
the exponential model is invariant under this extension. On the other hand, if $\mathcal{X}$ is the logistic model,
then $\mathcal{Y}_1$ is the type II generalized logistic model, which will also be invariant under this extension.

In general, when a statistical model is invariant under an extension that is stable,
it is because that model can be obtained as the extension of a submodel of it.

The extension of $\mathcal{X}$ obtained through minima and precursors of minima is:
$$  \mathcal{Y}_2 = min_{\mathcal{N}}(min^{-1}_{\mathcal{N}}(\mathcal{X})) =
\{Y_{\theta,\eta}: F_{Y_{\theta,\eta}} = (F_{X_{\theta}})^{e^{-\eta}}, ~ \theta \in \Theta, \eta \in (-\infty,\infty) \}, $$
which is a special case of the extension in Cordeiro et al (2013).
When one restricts $\eta \ge 0$, one obtains
$\mathcal{Y}_2^{\prime} = min_{\mathcal{N}}(\mathcal{X})$, and when one restricts $\eta \le 0$ one obtains
$\mathcal{Y}_2^{\prime \prime} = min^{-1}_{\mathcal{N}}(\mathcal{X})$, and therefore here
$\mathcal{Y}_2 = min_{\mathcal{N}}(\mathcal{X}) \cup min^{-1}_{\mathcal{N}}(\mathcal{X})$.

In this case, if $\mathcal{X}$ is for example the Gumbel model with the location parameter fixed, then $\mathcal{Y}_2$ is the
two parameter Gumbel model, and because of the stability of this extension the two parameter Gumbel model
model will be invariant under this extension. On the other hand,
when $\mathcal{X}$ is the logistic model then $\mathcal{Y}_2$ becomes the type I generalized logistic model which
by stability will also be invariant under this extension.

\textbf{Example 8.2:} Let the stopping model be the zero-truncated geometric in Example 6.2,
$$  \mathcal{N} = \{ N_{\eta}: h_{N_{\eta}} = \frac{t}{(1-t)e^{\eta} + t}, ~ \eta \in [0,\infty) \},  $$
which is extreme auto-reversible and includes $N_I$.

As a consequence of this auto-reversibility the extension of $\mathcal{X}$ obtained through maxima and their precursors or
through minima and their precursors here coincide, and it is:
$$  \mathcal{Y} =  \mathcal{Y}_1 = \mathcal{Y}_2 =
\{Y_{\theta,\eta}: F_{Y_{\theta,\eta}} =  \frac{F_{X_{\theta}}}{(1-F_{X_{\theta}})e^{\eta}+F_{X_{\theta}}}, ~ \theta \in \Theta, \eta \in (-\infty,\infty) \},  $$
which is the Marshall-Olkin extension of $\mathcal{X}$.
There is a huge literature using this model extension.
Here, when one restricts $\eta \ge 0$ one obtains $\mathcal{Y}^{\prime} = max_{\mathcal{N}}(\mathcal{X}) = min^{-1}_{\mathcal{N}}(\mathcal{X})$,
and when one restricts $\eta \le 0$ one obtains
$\mathcal{Y}^{\prime \prime} = min_{\mathcal{N}}(\mathcal{X}) = max^{-1}_{\mathcal{N}}(\mathcal{X})$.
As a consequence, this is the only example considered here where $\mathcal{Y} = max_{\mathcal{N}}(\mathcal{X}) \cup min_{\mathcal{N}}(\mathcal{X})$.

When for example $\mathcal{X}$ is the logistic model with the location parameter fixed
the extended model, $\mathcal{Y}$, is the two parameter logistic model.
Because of statistical stability of this extension, applying it again, now on the two-parameter logistic model, leaves the
model unchanged, which means that this two parameter model is invariant under this extension.

\textbf{Example 8.3:} Let the stopping model be the $\mathcal{N}_{\alpha}$ in Example 6.3 for a given $\alpha \ge 0$.
Like the geometric model, this one is also extreme auto-reversible, but it only includes $N_I$ when $\alpha=0$,
which is when it becomes the geometric model.

As a consequence of this auto-reversibility, the extensions of $\mathcal{X}$ obtained through maxima and their precursors and the
ones obtained through minima and their precursors coincide and are:
$$  \mathcal{Y}_{\alpha} = \{Y_{\theta,\eta}: F_{Y_{\theta,\eta}} = \frac{1}{\alpha}\ln\left(1+\frac{(e^{\alpha F_{X_{\theta}}}-1)(e^{\alpha}-1)}{(e^{\eta}-1)(e^{\alpha}-e^{\alpha F_{X_{\theta}}})+e^{\alpha}-1}\right), ~ \theta \in \Theta, \eta \in (-\infty,\infty) \}, $$
with $\mathcal{Y}^{\prime}_{\alpha} = max_{\mathcal{N}_{\alpha}}(\mathcal{X}) = min^{-1}_{\mathcal{N}_{\alpha}}(\mathcal{X})$
when one restricts $\eta \ge \alpha$, and with
$\mathcal{Y}^{\prime \prime}_{\alpha} = min_{\mathcal{N}_{\alpha}}(\mathcal{X}) = max^{-1}_{\mathcal{N}_{\alpha}}(\mathcal{X})$
when one restricts $\eta \le - \alpha$.
When one restricts $\eta \in (-\alpha,\alpha)$ one obtains $\mathcal{Y}^{\prime \prime \prime}_{\alpha} =
max_{\mathcal{N}_{\alpha}}(max^{-1}_{\mathcal{N}_{\alpha}}(\mathcal{X}))= max_{\mathcal{N}_{\alpha}}(min_{\mathcal{N}_{\alpha}}(\mathcal{X}))$
with $\eta_1,\eta_2$ such that $\eta_2 - \eta_1 \in (-\alpha,\alpha)$,
but this restricted transformation does not coincide with any of the transformations in Definition 2.

Different from what happens under the geometric model with $\alpha=0$, when $\alpha>0$
neither $max_{\mathcal{N}_{\alpha}}(\mathcal{X})$ nor $min_{\mathcal{N}_{\alpha}}(\mathcal{X})$
work as a model extension of $\mathcal{X}$, and
$max_{\mathcal{N}_{\alpha}}(\mathcal{X})\cup min_{\mathcal{N}_{\alpha}}(\mathcal{X}) \subset \mathcal{Y}_{\alpha}$
with an inclusion often strict. Hence this is an example where
$max_{\mathcal{N}_{\alpha}}(\mathcal{X})\cup min_{\mathcal{N}_{\alpha}}(\mathcal{X})$ does
not work as a model extension of $\mathcal{X}$, but where using the $\mathcal{Y}_{\alpha}$ from Definition 3 does.

\textbf{Example 8.4:} Let the stopping model be the $\mathcal{N}_{\alpha,\beta}$ in Example 6.4 for a given $\alpha \ge 0$
and $\beta \ge 1$. Here $N_I$ is not in the model, and the model
is not auto-reversible and therefore it yields two different model extension mechanisms.

The extension of $\mathcal{X}$ obtained through maxima and their precursors is:
$$  \mathcal{Y}_{1_{\alpha,\beta}} =  \{Y_{\theta,\eta}: F_{Y_{\theta,\eta}} = \frac{1-\left(\frac{\left(1-{\rm e}^{\alpha+\eta}\right)\left(1-F_{X_{\theta}}\left(1-{\rm e}^{-\frac{\alpha}{\beta}}\right)\right)^{\beta}+{\rm e}^{\eta}-1}{\left({\rm e}^{\alpha}-{\rm e}^{\alpha+\eta}\right)\left(1-F_{X_{\theta}}\left(1-{\rm e}^{-\frac{\alpha}{\beta}}\right)\right)^{\beta}+{\rm e}^{\eta}-{\rm e}^{\alpha}}\right)^{\frac{1}{\beta}}}{1-{\rm e}^{-\frac{\alpha}{\beta}}},
~ \theta \in \Theta, \eta \in (-\infty,\infty) \},  $$
and here one obtains $\mathcal{Y}_{1_{\alpha,\beta}}^{\prime} = max_{\mathcal{N}_{\alpha,\beta}}(\mathcal{X})$ when $\eta \ge \alpha$, and
$\mathcal{Y}_{1_{\alpha,\beta}}^{\prime \prime} = max^{-1}_{\mathcal{N}_{\alpha,\beta}}(\mathcal{X})$ when $\eta \le - \alpha$.
When $\eta \in (-\alpha,\alpha)$ one has that
$\mathcal{Y}_{1_{\alpha,\beta}}^{\prime \prime \prime} = max_{\mathcal{N}_{\alpha,\beta}}(max^{-1}_{\mathcal{N}_{\alpha,\beta}}(\mathcal{X}))$
with $\eta_1,\eta_2$ such that $\eta_2 - \eta_1 \in (-\alpha,\alpha)$,
which can neither be obtained through $max_{\mathcal{N}_{\alpha,\beta}}(\cdot)$ nor through $max^{-1}_{\mathcal{N}_{\alpha,\beta}}(\cdot)$.

The model extension of $\mathcal{X}$ obtained through minima and their precursors is:
$$  \mathcal{Y}_{2_{\alpha,\beta}} = \{Y_{\theta,\eta}: F_{Y_{\theta,\eta}} = \frac{1-\left(\frac{\left({\rm e}^{\eta}-{\rm e}^{\alpha}\right)\left(1+F_{X_{\theta}}\left({\rm e}^{\frac{\alpha}{\beta}}-1\right)\right)^{\beta}+{\rm e}^{\alpha}-{\rm e}^{\alpha+\eta}}{\left({\rm e}^{\eta}-1\right)\left(1+F_{X_{\theta}}\left({\rm e}^{\frac{\alpha}{\beta}}-1\right)\right)^{\beta}-{\rm e}^{\eta+\alpha}+1}\right)^{\frac{1}{\beta}}}{1-{\rm e}^{\frac{\alpha}{\beta}}}, ~ \theta \in \Theta, \eta \in (-\infty,\infty) \}, $$
and one obtains $\mathcal{Y}_{2_{\alpha,\beta}}^{\prime} = min_{\mathcal{N}_{\alpha,\beta}}(\mathcal{X})$ when $\eta \ge \alpha$, and
$\mathcal{Y}_{2_{\alpha,\beta}}^{\prime \prime} = min^{-1}_{\mathcal{N}_{\alpha,\beta}}(\mathcal{X})$ when $\eta \le - \alpha$.
When $\eta \in (-\alpha,\alpha)$ one has that
$\mathcal{Y}_{2_{\alpha,\beta}}^{\prime \prime \prime} = min_{\mathcal{N}_{\alpha,\beta}}(min^{-1}_{\mathcal{N}_{\alpha,\beta}}(\mathcal{X}))$
with $\eta_1,\eta_2$ such that $\eta_2 - \eta_1 \in (-\alpha,\alpha)$,
which can neither be obtained through $min_{\mathcal{N}_{\alpha,\beta}}(\cdot)$ nor through $min^{-1}_{\mathcal{N}_{\alpha,\beta}}(\cdot)$.

Here
$max_{\mathcal{N}_{\alpha,\beta}}(\mathcal{X})\cup max^{-1}_{\mathcal{N}_{\alpha,\beta}}(\mathcal{X}) \subset \mathcal{Y}_{1\alpha,\beta}$,
and $min_{\mathcal{N}_{\alpha,\beta}}(\mathcal{X})\cup min^{-1}_{\mathcal{N}_{\alpha,\beta}}(\mathcal{X}) \subset \mathcal{Y}_{2\alpha,\beta}$
with these inclusions being most often strict.

\textbf{Example 8.5:} Let the stopping model be the $\mathcal{N}_{\alpha,n}$ in Example 6.5 for a given $\alpha \ge 0$
and $n \in \mathbb{N}^{+}$. This model does not include $N_I$ and it is not extreme auto-reversible, but it is reversible with
the $\mathcal{N}_{\alpha,\beta=n}$ in Example 6.4. As a consequence, the model extension
obtained with $\mathcal{N}_{\alpha,n}$ through maxima and precursors of maxima,
coincide with the model extension obtained
with the $\mathcal{N}_{\alpha,\beta=n}$ of Example 6.4 through minima and precursors of minima, and viceversa.

\section{Final comments}

The main contribution of this article is putting together a set of new concepts
needed to define and untangle the properties of a large family of statistical model transformation mechanisms
that lead to statistical models useful for the analysis of extreme value data and in reliability.
The concepts introduced are:
\begin{enumerate}\vspace{-.3cm}
\item{} the notion of N-extreme precursors, which can be understood as the inverse
of $N$-stopped maxima and minima, and the model extension mechanisms derived from them (Definitions 1, 2 and 3),
which help generalize Marshall-Olkin extensions beyond geometric stopping,
\item{} the concept of statistical stability of a statistical model extension (Definition 5),
which applies to any statistical model extension and not
just to the ones considered in this paper,
\item{}
the idea of extreme reversible and auto-reversible stopping models (Definitions 6 and 7),
under which the extensions based on randomly stopped maxima and their inverses coincide with the
extensions based on randomly stopped minima and their inverses,
\item{}
and the idea of stopping models closed under pgf composition (Definition 8), which are the ones leading to
statistically stable randomly stopped extreme type of extensions.
\end{enumerate}\vspace{-.3cm}
All these new concepts are needed for the picture to be complete. In particular, if we
touch on methods to generate stopping models that are auto-reversible and/or closed under pgf composition
other than the geometric model, it is to help understand that the role played by geometric stopping is not
as unique as one might think after reading Marshall Olkin (1997).

A second contribution of this article are a set of theoretical results stating that uniparametric stopping models closed under pgf
composition can always be parametrized through $\theta=\Pr(N=1)$ with a parameter space
of the form $(0,\theta_0]$
(Theorem 1), and that the pgfs of these models commute under composition among themselves and with their inverses (Theorems 2 and 3).
These results are then used in Section 7 to determine conditions leading to statistically stable extensions.

Only two of the families of statistically stable model extensions presented in Section 8
are based on stopping models that are both closed
under pgf composition and extreme auto-reversible.
And the geometric model is the only stopping model that we know that shares these two features and includes $N_I$.
Nevertheless, note that in order to obtain statistically stable extensions through Definition 3,
one only needs that the stopping model be closed under pgf composition.

The only consequence of using stopping models that, unlike the geometric model,
are not extreme auto-reversible is that the extension based on maxima and
their inverses does not coincide with the extension based on minima and their inverses,
and using stopping models that, unlike geometric, do not include $N_I$
does neither affect the statistical stability nor the fact that the transformations presented
in Definition 3 always work as an extension.

Finally, note that our definition of statistical stability is
extremely basic and fundamental. A statistical model transformation
is statistically stable only if using that transformation twice in a row on any statistical model has the
same effect as using that transformation just once. The only reason that we can think for not finding the notion
of statistical stability anywhere in the statistical literature is that it might be difficult to prove results of that kind outside the
specific context of randomly stopped extreme transformations, and the closely related area of
randomly stopped sum transformations; It is easy to check that stopping models closed under pgf composition also
lead to randomly stopped sum model extensions that are statistically stable.

We consider statistical stability to be a property that should be central in the study
of any type of statistical model extension and not just in the study of the specific
extensions considered here, and we intend to keep investigating that.

\section*{Appendix 1: Model extensions when $\mathcal{N}$ is the zero truncated Poisson or the logarithmic model}

The zero-truncated Poisson($\alpha$) model is defined through the set of pgfs:
$$ \mathcal{N} = \{ N_{\alpha}: h_{N_{\alpha}} = \frac{e^{\alpha t} - 1}{e^{\alpha} - 1}, ~~ \alpha \in [0,\infty) \}. $$
This model includes $N_I$
and therefore both the basic transformations in Definition 2 as well as the combined transformations
in Definition 3 are extensions,
but this model is neither extreme auto-reversible, because $\Pr[N_{\alpha}=1] \ne 1/E[N_{\alpha}]$, nor closed under pgf composition, because
$$ \Pr( N_{\alpha} = 1) =  \frac{\alpha}{e^{\alpha}-1}, $$
$$ \Pr( N_{\alpha} = 2) =  \frac{1}{2}\frac{\alpha^{2}}{e^{\alpha}-1}, $$
and therefore it does not satisfy the necessary condition of Corollary 2 for being closed,
$$   \frac{\Pr(N_{\alpha}=2)}{\Pr(N_{\alpha}=1)(1-\Pr(N_{\alpha}=1))} =
\frac{\alpha}{2}+\frac{1}{2}\frac{\alpha^{2}}{e^{\alpha}-\alpha-1} \ne Constant. $$
The four basic
extensions of $\mathcal{X} = \{X_{\theta}: F_{X_{\theta}}, ~ \theta \in \Theta \}$ obtained through Definition 2 are,
$$  max_{\mathcal{N}}(\mathcal{X}) = \{Y_{\theta,\alpha}: F_{Y_{\theta,\alpha}} =  \frac{e^{\alpha F_{X_{\theta}}} - 1}{e^{\alpha} - 1}, ~ \theta \in \Theta, \alpha \in [0,\infty) \},   $$
$$  max^{-1}_{\mathcal{N}}(\mathcal{X}) = \{Y_{\theta,\alpha}: F_{Y_{\theta,\alpha}} = \frac{\ln(1+(e^{\alpha}-1)F_{X_{\theta}})}{\alpha}, ~ \theta \in \Theta, \alpha \in [0,\infty) \},  $$
$$  min_{\mathcal{N}}(\mathcal{X}) = \{Y_{\theta,\alpha}: F_{Y_{\theta,\alpha}} = \frac{e^{\alpha}(1 - e^{-\alpha F_{X_{\theta}}})}{e^{\alpha} - 1}, ~ \theta \in \Theta, \alpha \in [0,\infty) \},  $$
$$  min^{-1}_{\mathcal{N}}(\mathcal{X}) = \{Y_{\theta,\alpha}: F_{Y_{\theta,\alpha}} = 1 - \frac{\ln(1+(e^{\alpha}-1)(1-F_{X_{\theta}}))}{\alpha}, ~ \theta \in \Theta, \alpha \in [0,\infty) \},  $$
and the four combined extensions of $\mathcal{X}$ obtained through Definition 3 are,
$$  max_{\mathcal{N}}(max^{-1}_{\mathcal{N}}(\mathcal{X})) = $$
$$ \{Y_{\theta,\alpha_1,\alpha_2}: F_{Y_{\theta,\alpha_1,\alpha_2}} = \frac{1}{{\rm e}^{\alpha_2}-1}\left(\left(1+\left({\rm e}^{\alpha_1}-1\right)F_{X_{\theta}}\right)^{\frac{\alpha_2}{\alpha_1}}-1\right), ~ \theta \in \Theta, \alpha_1,\alpha_2 \in [0,\infty)\},   $$
$$  max^{-1}_{\mathcal{N}}(max_{\mathcal{N}}(\mathcal{X})) = $$
$$ \{Y_{\theta,\alpha_1,\alpha_2}: F_{Y_{\theta,\alpha_1,\alpha_2}} = \frac{1}{\alpha_2}\ln\left(1+\frac{\left({\rm e}^{\alpha_2}-1\right)\left({\rm e}^{\alpha_1\,F_{X_{\theta}}}-1\right)}{{\rm e}^{\alpha_1}-1}\right), ~ \theta \in \Theta, \alpha_1,\alpha_2 \in [0,\infty)\},  $$
$$  min_{\mathcal{N}}(min^{-1}_{\mathcal{N}}(\mathcal{X})) = $$
$$ \{Y_{\theta,\alpha_1,\alpha_2}: F_{Y_{\theta,\alpha_1,\alpha_2}} = \frac{{\rm e}^{\alpha_2}}{{\rm e}^{\alpha_2}-1}\left(1-\left(1-\left(1-{\rm e}^{-\alpha_1}\right)F_{X_{\theta}} \right)^{\frac{\alpha_2}{\alpha_1}}\right), ~ \theta \in \Theta, \alpha_1,\alpha_2 \in [0,\infty)\},  $$
$$  min^{-1}_{\mathcal{N}}(min_{\mathcal{N}}(\mathcal{X})) = $$
$$ \{Y_{\theta,\alpha_1,\alpha_2}: F_{Y_{\theta,\alpha_1,\alpha_2}} = 1-\frac{1}{\alpha_2}\ln\left(\frac{\left({\rm e}^{\alpha_2}-1\right)\left({\rm e}^{\alpha_1\,\left(1-F_{X_{\theta}}\right)}-1\right)}{{\rm e}^{\alpha_1}-1}+1\right), ~ \theta \in \Theta, \alpha_1,\alpha_2 \in [0,\infty)\}.  $$
Furthermore, according to Example 5.2 the Logarithmic$(p)$ model defined through:
$$  \mathcal{N} =  \{ N_{p}: h_{N_{p}} = \frac{\log{(1-pt)}}{\log{(1-p)}}, ~~ p \in [0,1) \},  $$
where $p = 1 - e^{-\alpha}$, is extreme reversible with the zero-truncated Poisson model.
As a consequence of that property, the set of model extensions
obtained through Definitions 2 and 3 using the Logarithmic$(p)$ model
coincide with the set of extensions
obtained using the zero-truncated Poisson model presented in this Appendix.

The specific extensions for the Logarithmic$(p)$
model are the ones listed above for the truncated Poisson model after replacing $\alpha$ by $-log(1-p)$,
and after switching $max_{\mathcal{N}}$ and $min^{-1}_{\mathcal{N}}$
and switching $min_{\mathcal{N}}$ and $max^{-1}_{\mathcal{N}}$. For example the $max_{\mathcal{N}}(\cdot)$ transformation when
$\mathcal{N}$ is Logaritmic$(p)$ is the $min^{-1}_{\mathcal{N}^{'}}(\cdot)$ transformation when $\mathcal{N}^{'}$
is truncated Poisson$(\alpha=-log(1-p))$, the $min_{\mathcal{N}}(\cdot)$ transformation is the $max^{-1}_{\mathcal{N}^{'}}(\cdot)$
transformation, and so on.

\section*{Appendix 2: Proof of Theorem 1 in Section 6}

Let $\mathcal{N}
= \{N_{\delta}: h_{N_{\delta}}= \sum_{i=1}^{\infty} p_i(\delta) t^i, ~ \delta \in \mathcal{D} \}$, with $p_i(\delta)=\Pr(N_{\delta}=i)$,
be a uniparametric stopping model closed under pgf composition and with a continuously differentiable parametrization and
a connected parameter space with a non-empty interior.
That is, assume that the stopping model, besides being closed under pgf composition is also such that $p_i(\delta)$ is
continuously differentiable in
$\delta$ for any $i$, and that $\mathcal{D}$ is a non-empty interval of $\mathbb{R}$, which is
what in the paper is denoted in short as a model ``uniparametric and closed under pgf composition."

Let $k_1 = \min \{k:\exists \delta \in \mathcal{D} ~ \mbox{such that} ~ p_k(\delta) > 0 \}$, and let
$\mathcal{D}_{k} = \{ \delta \in \mathcal{D}: ~ p_{k}(\delta)>0 ~ \mbox{and} ~ p_j(\delta) = 0 ~ \mbox{for all} ~ j < k \}$ and
$\mathcal{N}_k = \{ N_{\delta} \in \mathcal{N}: ~ \delta \in \mathcal{D}_{k} \}$
be partitions of $\mathcal{D}$ and $\mathcal{N}$.
Even though it is assumed that $\mathcal{N} = \cup_{k=1}^{\infty}\mathcal{N}_{k}$ is closed under pgf composition
and $\mathcal{D} = \cup_{k=1}^{\infty} \mathcal{D}_{k}$ is connected, the $\mathcal{N}_{k}$
do not have to be closed under pgf composition and the $\mathcal{D}_{k}$ do not have to be connected.

Let $N_{\tilde{\delta}} = N_{\delta_1} \circ N_{\delta_2}$ stand for the random variable in $\mathcal{N}$
with pgf $h_{N_{\tilde{\delta}}} = h_{N_{\delta_1}}\circ h_{N_{\delta_2}},$ and denote $\tilde{\delta} = \delta_1 \circ \delta_2$.
Note that if $N_{\delta_1} \in \mathcal{N}_{k_1}$
and $N_{\delta_2} \in \mathcal{N}_{k_2}$, then $N_{\tilde{\delta}} = N_{\delta_1} \circ N_{\delta_2} \in \mathcal{N}_{k_1 k_2}$.

Let $\mathcal{N}^{\circ m} = \{ N_{\delta_1} \circ N_{\delta_2} \circ \cdots \circ N_{\delta_m}: ~ N_{\delta_j} \in \mathcal{N} \}$,
where $m \in \mathbb{N}$, and let $N^{\circ m}_{\delta} = N_{\delta^{\circ m}} = N_{\delta} \circ N_{\delta} \circ \cdots \circ N_{\delta}$.
Note that if $\mathcal{N}$ is uniparametric and closed under pgf composition, then $\mathcal{N}^{\circ m} \subset \mathcal{N}$,
and $\mathcal{N}^{\circ m}$ is also a uniparametric and closed under pgf composition.

Next we state and prove four lemmas required to prove Theorem 1.

Stating that ``$\Pr(N_{\delta}=1)>0$ for all $N_{\delta} \in \mathcal{N}$", as in the first point of Theorem 1,
is equivalent to stating that ``$\mathcal{D}_1 = \mathcal{D}$", because if $k_1 > 1$, or if
$k_1=1$ but there exists a $\delta_0 \in \mathcal{D}$ such that $p_{k_1}(\delta_0)=0$, both statements are false,
and otherwise they are both true.

Theorem 1 will be proven by \emph{reductio ad absurdum}, and therefore one needs
to know how would the stopping models be if they were ``uniparametric and
closed under pgf composition, but $\Pr(N_{\delta}=1)=0$ for some $N_{\delta}\in\mathcal{N}$",
and hence if ``$\mathcal{D}_1 \subsetneq \mathcal{D}$".
That is the purpose of the first lemma.

\begin{lemma}
If $\mathcal{N}$ is ``uniparametric
and closed under pgf composition, but such that $\Pr(N_{\delta}=1)=0$
for some $N_{\delta} \in \mathcal{N}$", then
at the boundary of any connected component of $\mathcal{D}_{k_1}$, denoted $\mathcal{D}_{k_1}^{*}$,
there is a $\delta_0^{*}$, from $\mathcal{\mathring{D}}$, with $p_{k_1}(\delta_0^{*})=0$ and
such that if a sequence $\delta_i$ in $\mathcal{D}_{k_1}^{*}$ converges to $\delta_0^{*}$,
then $p_{k_1}(\delta_i)$ converges to $0$.
\end{lemma}
{\sl{Proof:}}
Under the conditions stated in the lemma, there exists a $\delta_0$ in $\mathcal{D}$ such that
$p_1(\delta_0) = 0$, and so such that if $\delta \in \mathcal{D}_{k_1}^{*}$ then
$N_{\tilde{\delta}} = N_{\delta_0} \circ N_{\delta}$ belongs to $\mathcal{N}_k$ with $k > k_1$.
As a consequence, $p_{k_1}(\tilde{\delta})=0$ and therefore $\mathcal{D}_{k_1} \subsetneq \mathcal{D},$
or what is equivalent, $\mathcal{D} \setminus \mathcal{D}_{k_1} \neq \emptyset$.
Given that for any $\delta$ in $\mathcal{D}\setminus\mathcal{D}_{k_1}$ one has that
$p_{k_1}(\delta)=0$, any connected component $\mathcal{D}_{k_1}^{*} \subsetneq \mathcal{D}$
will be an open interval with at least one extreme $\delta_{0}^{*}$ in $\mathring{\mathcal{D}} \setminus \mathcal{D}_{k_1}$,
and therefore with $p_{k_1}(\delta_{0}^{*})=0$, (if the two extremes of $\mathcal{D}_{k_1}^{*}$ were the ones of $\mathcal{D}$,
then $\mathcal{D}_{k_1}^{*} = \mathcal{D}$ but $\mathcal{D}_{k_1}^{*} \subsetneq \mathcal{D}$).
By the continuity of $p_{k_1}(\delta)$, if a sequence $\delta_i$
in $\mathcal{D}_{k_1}^{*}$ converges to $\delta_{0}^{*}$, then $p_{k_1}(\delta_i)$ converges to $0$.
%
%
\hfill{$\square$}

\begin{remark}
When ``$\Pr(N_{\delta}=1)>0$ for all $N_{\delta} \in \mathcal{N}$," as in the first
statement in Theorem 1, and therefore when
``$\mathcal{D}_1 = \mathcal{D} = \mathcal{D}_{k_1}$," one also has that
$\mathcal{D}_1$ is connected and that $\inf_{\delta \in \mathcal{D}_1} p_1(\delta) = 0$, because
$N_{\delta}^{\circ m} = N_{\delta^{\circ m}} \in \mathcal{N}$
for any $\delta \in \mathcal{D}_1$ and any $m \in \mathbb{N}$, and
$\Pr(N_{\delta}^{\circ m}) = p_1(\delta^{\circ m}) = p_1(\delta)^{m}$
which tends to $0$ when $m$ tends to infinity.
\end{remark}

From now on, $\mathcal{D}_{k_1}^{*}$ will stand for a connected component of $\mathcal{D}_{k_{1}}$,
that is, an open interval with $\delta_0^{*}$ as at least one of its endpoints.
Different from Lemma 1, the next three lemmas relate to properties of all models
``uniparametric and closed under pgf composition".
\begin{lemma}
If $\mathcal{N}$ is ``uniparametric and closed under pgf composition", then one can not have
$p_{k_1}(\delta) = C \neq 0$ for $\delta \in (\delta_{1}-\epsilon,\delta_{1}+\epsilon) \subset \mathcal{D}_{k_1}^{*}$
with $\epsilon>0$.
\end{lemma}
{\sl{Proof:}} The result will be proven by \emph{reductio ad absurdum}, checking that if $p_{k_1}(\delta) = C \neq 0$ in
an interval in $\mathcal{D}_{k_1}^{*}$, the model can not be uniparametric.

If $p_k(\delta) = C_k \ne 0$ and so $p_k^{\prime}(\delta)=0$
for $\delta \in (\delta_{1}-\epsilon,\delta_{1}+\epsilon)$ and for all $k > k_1$,
the parametrization $\delta$ would not be identifiable.
Let $m > k_1$ be the first term such that
there exists a $\delta^{\prime}_1 \in (\delta_1 - \epsilon, \delta_1 + \epsilon)$ with
$p_m^{\prime}(\delta^{\prime}_1) \neq 0$ and so by the continuity of the derivative such that $p_m(\delta)$ is strictly monotone in an environment
$(\delta^{\prime}_1 - \epsilon^{\prime}, \delta^{\prime}_1 + \epsilon') \subset (\delta_1 - \epsilon, \delta_1 + \epsilon)$,
and let $\delta_1^{\prime}$ and $\epsilon^{\prime}$ be such that
$p_k(\delta) = C_k$  in $(\delta_1^{\prime} - \epsilon^{\prime}, \delta_1^{\prime} + \epsilon^{\prime})$
for all $k_1 \le k < m$. (If $p_m^{\prime}(\delta_1) \neq 0,$ one can choose $\delta^{\prime}_1 = \delta_1$).
This means that in this interval one can parametrize the $N_{\delta}$
through $\theta = p_m(\delta)$. From now on, we relabel $\delta_1^{\prime}$ by $\delta_1$, $\epsilon^{\prime}$ by $\epsilon$ and
the $C$ of the statement of the lemma by $C_{k_1}$.

If $d_1, d_2$ are in $(\delta_1 - \epsilon, \delta_1 + \epsilon)$, then $N_{d_1} \circ N_{d_2}$
belongs to $\mathcal{N}$ and its pgf is:
$$ h_{N_{d_1}}\circ h_{N_{d_2}}(t) = p_{k_1}(d_1)
(h_{N_{d_2}}(t))^{k_1} + \cdots + p_m(d_1)(h_{N_{d_2}}(t))^m + \cdots$$
The first term of this pgf where $\theta_2 = p_m(d_2)$ appears is
$$p_{k_1}(d_1)(h_{N_{d_2}}(t))^{k_1} = p_{k_1}(d_1)(p_{k_1}(d_2)t{}^{k_1} + \cdots + p_m(d_2)t{}^m + \cdots)^{k_1} = $$
$$ C_{k_1}(C_{k_1} t^{k_1} + \cdots + C_{m-1} t^{m-1} + \theta_2 t{}^m + \cdots + p_j(\theta_2)t^j + \cdots)^{k_1}, $$
and it is under the form
$$ k_1 C_{k_1}(C_{k_1}t^{k_1})^{k_1-1}\theta_2 t{}^m, $$
and therefore $\theta_2 = p_m(d_2)$ appears first in the coefficient of the term $t^{k_1^2+m-k_1}$ of the pgf of $N_{d_1} \circ N_{d_2}$,
with a coefficient equal to
$k_1 C^{k_1} \theta_2 + Q$,
where $Q$ is a constant term that depends on $p_k(d_1) = p_k(d_2) = C_k$ for $k_1 \le k < m$. Therefore that first term only depends on $\theta_2$.

Analogously, the first term of the pgf of $N_{d_1} \circ N_{d_2}$ where $\theta_1 = p_m(d_1)$ appears is
$p_{m}(d_{1})(h_{N_{d_{2}}}(t))^{m}$, and it is under the form
$p_m(d_1)(p_{k_1}(d_2)t^{k_1})^m$,
and therefore $\theta_1 = p_m(d_1)$ appears first in the coefficient of the term $t^{mk_1}$ of that pgf, with a coefficient equal to
$C^m \theta_1 + R(\theta_2)$, where $R(\theta_2)$ are terms that depend on
$p_k(d_1) = C_k$ for $k_1 \le k < m$ and on $\theta_2 = p_m(d_2)$.

Given that $\theta_1$ and $\theta_2$ are not related,
one needs two parameters to describe $N_{d_1} \circ N_{d_2} \in \mathcal{N}$, and $\mathcal{N}$ can not be uniparametric, which proves the result.
\hfill{$\square$}

\begin{remark}
One consequence of Lemma 2 is that $p_{k_1}^{\prime}(\delta)$ can not be zero in any interval, and
therefore the number of critical points of $p_{k_1}(\delta)$ in $\mathcal{D}_{k_1}^{*}$, with $p_{k_1}^{\prime}(\delta) = 0$, is either
finite or countable.
\end{remark}

\begin{lemma}
If $\mathcal{N}$ is ``uniparametric and closed under pgf composition", then
$p_{k_1}(\delta)$ is strictly monotone in $\mathcal{D}_{k_1}^{*}$, and
therefore the $N_{\delta}$ with $\delta \in \mathcal{D}_{k_1}^{*}$ can be parametrized through $\theta = p_{k_1}(\delta)$.
\end{lemma}
{\sl{Proof:}}
One assumes that $p_{k_1}(\delta)$ is continuously differentiable and so continuous, with
$0 < p_{k_1}(\delta)\le 1$ for all $\delta \in \mathcal{D}_{k_1}^*$.
By Lemma 2 one knows that $p_{k_1}(\delta)$ is not constant in
any subinterval of $\mathcal{D}_{k_1}^{*}$, and by Lemma 1 and Remark 1
one knows that one always has that $\inf_{\delta \in \mathcal{D}_{k_1}^{*}} p_1(\delta) = 0$.

To prove that $p_{k_1}(\delta)$ is strictly monotone in $\mathcal{D}_{k_1}^{*}$, one needs to prove that there are no local
maxima or local minima in its interior, and no fluctuating critical points in its closure, where
a fluctuating critical point $\delta_f$ is the limit of critical points such that $p_{k_1}(\delta)$ fluctuates when $\delta$ converges to $\delta_f$.
In our case though, it will be sufficient to
check that there is no local maximum of $p_{k_1}(\delta)$ in the interior of $\mathcal{D}_{k_1}^{*}$, because:
\begin{itemize}\vspace{-.3cm}
\item{}
if there is no local maximum there can not be any fluctuating critical points $\delta_f$,
because $\delta_f$ would have to be the accumulation point of fluctuating critical points
(the points where $p_{k_1}(\delta)$ fluctuates from increasing to decreasing
can neither be local maxima nor local minima), and as a consequence
the set of fluctuating critical points $\delta_f$ in $\bar{\mathcal{D}}_{k_1}^{*}$ would be a perfect set and the number of critical points in $\mathcal{D}_{k_1}^{*}$
would have a cardinality larger than countable, which is impossible by Remark 2,
\item{}
and if there is no local maximum in the interior of $\mathcal{D}_{k_1}^{*}$ and no fluctuating critical point in its closure,
there can not exist an isolated minimum in there
either, because if it did $p_{k_1}(\delta)$ would not be able to take values close enough to $0$ near the boundary of $\mathcal{D}_{k_1}^{*}$.
\end{itemize}\vspace{-.3cm}

That there is no local maximum is proven by \emph{reductio ad absurdum}.
Assume that $p_{k_1}(\delta)$ had a local maximum, $\delta_m$, with
$\delta_1 < \delta_m < \delta_2$ where $\delta_1, \delta_2 \in \mathcal{D}_{k_1}^{*}$, and by continuity of the derivative with
$p_{k_1}^{\prime}(\delta) \ge 0$ for $\delta \in (\delta_1 - \epsilon, \delta_m)$ and with $p_{k_1}^{\prime}(\delta) \le 0$
for $\delta \in (\delta_m, \delta_2 + \epsilon)$, and let
$q=p_{k_1}(\delta_1) = p_{k_1}(\delta_2) < p_{k_1}(\delta_m) = q_m$.
Because of the sign of the derivative, $p_{k_1}(\delta)$ would be increasing in
$(\delta_1-\epsilon,\delta_m)$ and by Lemma 2 it would be strictly increasing there, and
$p_{k_1}(\delta)$ would be strictly decreasing in $(\delta_m, \delta_2+\epsilon)$ for the same reasons. This would mean that, because of continuity,
this distributions can be parametrized through $p_{k_1} = p_{k_1}(\delta)$ both in $[\delta_1,\delta_m]$ as well as in $[\delta_m,\delta_2]$.

Let $\delta,\delta^{*} \in \mathcal{D}_{k_1}$, and consider the family of distributions
$N_{\tilde{\delta}} = N_{\delta}\circ N_{\delta^{*}} \in \mathcal{N}_{k_1} \circ \mathcal{N}_{k_1}\subset \mathcal{N}_{k_1^2}
\subset  \mathcal{N}^{\circ 2} \subset \mathcal{N}$,
which is also uniparametric, closed under pgf composition, and with a $k_1^{'}$ equal to $k_1^2$, and lets denote
$\tilde{\delta}(\delta,\delta^{*})$ by $\delta \circ \delta^{*}$.

If $\delta,\delta^{*} \in \mathcal{D}_{k_1}^{*}$, then  $\delta \circ \delta^{*} \in \mathcal{D}_{k_1^2}^{*}$, which is a connected subset of
$\mathbb{R}$, and if $\delta,\delta^{*} \in [\delta_1,\delta_2] \subset \mathcal{D}_{k_1}^{*}$, then
$\delta\circ\delta^{*} \in [\delta_1,\delta_2]^{\circ2} \subset \mathcal{D}_{k_1^2}^{*}$.
If $\delta \in [\delta_1, \delta_2]$, then $\{ \delta\circ\delta\} =
[\delta_1^{\circ2},\delta_2^{\circ2}]$
is a connected subset of $[\delta_1,\delta_2]^{\circ 2}$, because one knows that if
$\delta \neq \delta^{*}$ then $\delta^{\circ2}\neq\delta^{*\circ2}$, and that by continuity of the parametrization, if
$\delta_j \in (\delta_1, \delta_2)$ converges to $\delta$ in $[\delta_1,\delta_2]$, then
$\delta_j^{\circ 2}$ converges to $\delta^{\circ2}$. The set
$[\delta_1^{\circ2},\delta_2^{\circ2}]$, (or the set $[\delta_2^{\circ2},\delta_1^{\circ2}]$ if $\delta_2^{\circ2} < \delta_1^{\circ2}$),
can be parametrized through $\delta^{\circ2}$, and with that parametrization one has that
$p_{k_1^2}(\delta^{\circ 2}) = p_{k_1}(\delta) p_{k_1}(\delta)^{k_1} = p_{k_1}(\delta)^{k_1+1} \in [q^{k_1+1},q_m^{k_1+1}] \subset [0,1]$.
In particular, $p_{k_1^2}(\delta_1^{\circ 2}) = p_{k_1^2}(\delta_2^{\circ 2}) = q^{k_1+1}$ and
$p_{k_1^2}(\delta_m^{\circ 2}) = q_m^{k_1+1}$.

Consider now the set $\{ \delta_1 \circ [\delta_1,\delta_2] \} \subset [\delta_1,\delta_2]^{\circ 2}$, with a first
value $\delta_1^{\circ 2}$ that belongs to $[\delta_1^{\circ 2}, \delta_m^{\circ 2}]$.
Given that the set $\{ \delta_1 \circ [\delta_1,\delta_2] \}$ is connected and of dimension one, it has to be
contained in $[\delta_1^{\circ 2},\delta_m^{\circ 2}]$, because otherwise it would either contain $\delta_m^{\circ 2}$ or take
values smaller than $\delta_1^{\circ2}$, but
$p_{k_1^2}(\delta_1 \circ[\delta_1,\delta_2]) \subset [p_{k_1^2}(\delta_{\text{1}}^{\circ2}), p_{k_1^2}(\delta_m^{\circ 2}))$
and so its last value, $\delta_1 \circ \delta_2$, is in $[\delta_1^{\circ 2},\delta_m^{\circ 2}]$;
Given that $p_{k_1^2}(\delta_1 \circ \delta_2) = q^{k_1+1}$
and that $\delta_1^{\circ 2}$ is the unique point in $[\delta_1^{\circ 2},\delta_m^{\circ 2}]$
such that $p_{k_1^2}(\tilde{\delta}) = p_{k_1^2}(\delta_1 \circ \delta_2)=q^{k_1+1}$,
then $\delta_1 \circ \delta_2 = \delta_1 \circ \delta_1$. But this implies that
$h_{N_{\delta_1}}(h_{N_{\delta_2}}(t)) = h_{N_{\delta_1}}(h_{N_{\delta_1}}(t))$, and so that
$\delta_2 = \delta_1$, which is in contradiction with $\delta_1 < \delta_m < \delta_2$ and therefore there
can not exist local maxima of $p_{k_1}(\delta)$ in the interior of $\mathcal{D}_{k_1}^{*}$.

As a consequence, $\theta = p_{k_1}(\delta)$ is strictly monotone in $\mathcal{D}_{k_1}^{*}.$
%
\hfill{$\square$}

\begin{lemma}
If $\mathcal{N}$ is ``uniparametric and closed under pgf composition", then $k_1 = 1$.
\end{lemma}
{\sl{Proof:}}
We prove that if $k_1 > 1$, then $\mathcal{N}$ can not be ``uniparametric and closed under pgf composition"
by \emph{reductio ad absurdum}, in two steps.

In the first step we show that if $\mathcal{N}$ is ``uniparametric and closed under pgf composition" and $k_1 > 1$, then
$p_{k_1+j}(\theta) = \theta P_{j}(\theta)$ for all $j \ge 0$ and for every $\delta \in \mathcal{D}_{k_1}^*$,
where $\theta = p_{k_1}(\delta)$ and $P_j(\theta)$ is a polynomial of $\theta$.

In the second step we show that if $k_1 > 1$ and $\delta_i$ converges to $\delta_0$ and so $p_{k_1}(\delta_i)$ converges to $0$,
then $p_{k_1+j}(\delta_i)$ converges to $0$ for all $j \ge 0$, which given that $p_{r}(\delta_i)=0$ for all $r < k_1$ implies that
all the coefficients of $h_{\delta_0}(t)$ would be $0$
and so $\delta_0 \notin \mathcal{D}$, but that would contradict Lemma 1.
\begin{enumerate}\vspace{-.3cm}
\item{Step 1:}
We prove that if $\mathcal{N}$ is ``uniparametric and closed under pgf composition" and $k_1 > 1$,
then $p_{k_1+j}(\theta) = \theta P_{j}(\theta)$,
by induction.

It holds for $j=0$, because $p_{k_1}(\theta)=\theta P_0(\theta)$ with $P_0(\theta)=1$.
Lets assume that $p_{k_1+j}(\theta)=\theta P_j(\theta)$ for $j \le m-1$, and check what happens for $j=m$.

For any $\delta_1,\delta_2 \in \mathcal{D}_{k_1}^*$ one has that
$N_{\delta_1}\circ N_{\delta_2} \in \mathcal{N}^{\circ 2} \subset \mathcal{N}_{k_1^2} \subset \mathcal{N}$, and
$$h_{N_{\delta_1}\circ N_{\delta_2}}
=\sum_{r=0}^{\infty}\left(p_{k_1+r}(\delta_1)\left(\sum_{j=0}^{\infty}p_{k_1+j}(\delta_2)t^{k_1+j}\right)^{k_1+r}\right)$$
%
which by Lemma 3 can be parametritzed through $\tau=p_{k_1^2}(\delta)=p_{k_1}(\delta_1)p_{k_1}(\delta_2)^{k_1}=\theta_1\theta_2^{k_1}$,
and thus with $\theta_1 = \tau/\theta_2^{k_1}$.
The first term in $h_{N_{\theta_1}\circ N_{\theta_2}}$ where $p_{k_1+m}(\theta_2)$ appears is:
$$ p_{k_1}(\delta_1)k_1 (p_{k_1}(\delta_2)t^{k_1})^{k_1-1} p_{k_1+m}(\delta_2)t^{k_1+m}, $$
which is the term of order $k_1^2 + m$, and its coefficient is of the form
$$  p_{k_1^2 + m}(\tau) = \sum_r p_{k_1+r}(\theta_1)\prod_{i=1}^{k_1+r}p_{k_1+j_i}(\theta_2),  $$
where the sum is over all $r$ and $j_i$ such that $0 \le r$, $0 \le j_i$ and $0 \le \sum_{i=1}^{k_1+r} j_i = m - k_1 r$,
and therefore it can be written as:
$$ p_{k_1^2 + m}(\tau) = k_1 p_{k_1}(\theta_1)(p_{k_1}(\theta_2))^{k_1-1} p_{k_1+m}(\theta_2) + A, $$
where the first summand is the sum of the $k_1$ terms where $r=0$ and all $j_i=0$ except one $j_i$ that is equal to $m$,
and where $A$ is the sum of the remaining terms, all with $j_i < m$.
The first term where $p_{k_1+m}(\theta_1)$ appears is the one of order $k_1^2 + k_1 m$, larger than $k_1^2 + m$ because $k_1>1$,
and therefore the $p_{k_1+j}(\theta_1)$ with $j \ge m$ do not appear in the above expression for $p_{k_1^2 + m}(\tau)$.

Substituting $\theta_1 = \tau/\theta_2^{k_1}$ and using the fact that $p_{k_1+j}(\theta)=\theta P_j(\theta)$ when $j < m$,
the first sumand in the last expression for $p_{k_1^2 + m}(\tau)$ becomes:
$$ k_1 \theta_1 \theta_2^{k_1-1} p_{k_1+m}(\theta_2) =
k_1 \frac{\tau}{\theta_2} p_{k_1+m}(\theta_2),                $$
and the second sumand in that expression becomes:
$$ A = \sum_r p_{k_1+r}(\theta_1) \prod_{i=1}^{k_1+r} p_{k_1+j_i}(\theta_2) =
\sum_r \theta_1 P_r(\theta_1) \prod_{i=1}^{k_1+r} \theta_2 P_{j_i}(\theta_2) =  $$
$$  \sum_r \tau P_r \left(\frac{\tau}{\theta_2^{k_1}}\right) \frac{\theta_2^{r+k_1}}{\theta_2^{k_1}} \prod_{i=1}^{k_1+r} P_{j_i}(\theta_2),   $$
where the sum is over the set of $r$ and $j_i$ described above but excluding $j_i = m$.
Hence, the coefficient of the $k_1^2 + m$ term of $h_{N_{\delta_1}\circ N_{\delta_2}}$ is a polynomial
in $\tau$ of the form:
$$p_{k_1^2+m}(\tau)=\tau \left(\frac{k_1 p_{k_1+m}(\theta_2)}{\theta_2} + P_{m(1)}(\theta_2) \right) +
\sum_{i=2}^{c} \tau^i P_{m(i)}(\theta_2). $$
If $\mathcal{N}$ is uniparametric and closed under pgf composition, then
all the coefficients of $h_{N_{\delta_1}\circ N_{\delta_2}}$ have to depend only on $\tau$, and hence all
the coefficients of this polynomial in $\tau$ need to be constant.
In particular, for the coefficient of the first term in $\tau$ to be constant, equal to $C_m$, one needs that
$k_1 p_{k_1+m}(\theta_2) =
\theta_2(C_m - P_{m(1)}(\theta_2))$, and therefore that
$p_{k_1+m}(\theta) = \theta P_m(\theta)$.

\item{Step 2:}
If $k_1 > 1$ and $\theta_i = p_{k_1}(\delta_i)$ converges to $0$ when $\delta_i$ converges to $\delta_0$, then
it follows that $p_{k_1+m}(\delta_i) = p_{k_1}(\delta_i)P_m(p_{k_1}(\delta_i))$
will converge to $0$ for all $m \ge 0$,
and by definition of $k_1$, $p_j(\delta_i)=0$ for all $j<k_1$. Therefore $h_{\delta_0}(t) = 0$, which is not a pgf
and so $\delta_0 \notin \mathcal{D}$.
By Lemma 1, if $k_1 > 1$ there is one such $\delta_0 \in \mathring{\mathcal{D}}$ on the boundary of $\mathcal{D}_{k_1}^{*}$,
but we just proved that $\delta_0 \notin \mathring{\mathcal{D}}$, and therefore $k_1$ can not be larger than 1.
\hfill{$\square$}
\end{enumerate}\vspace{-.3cm}

{\sl{Proof of Theorem 1, in Section 6:}}

The fact that when $\mathcal{N}$ is ``\emph{uniparametric and closed under pgf composition}", then
$p_1(\delta) = \Pr(N_{\delta} = 1) > 0$ for all $\delta \in \mathcal{D}$,
follows from the fact that $k_1=1$ by Lemma 4, and that
if there existed a $\delta_0 \in \mathring{\mathcal{D}}$ with $p_1(\delta_0) = 0$, then the model $\mathcal{N}_0 =
\{ N_{\delta} \circ N_{\delta_0}: N_{\delta} \in \mathcal{N}\}$ would also be ``uniparametric and closed under pgf composition",
because $(N_{\delta_1}\circ N_{\delta_0})\circ(N_{\delta_2}\circ N_{\delta_0})=(N_{\delta_1}\circ N_{\delta_0}\circ N_{\delta_2})\circ N_{\delta_0}
\in\mathcal{N}_0$, but $\mathcal{N}_0$ would have $k_1 > 1$ because
$P(N_{\delta} \circ N_{\delta_0} = 1) = p_1(\delta_N) p_1(\delta_0) = 0$
for any $N_{\delta} \in \mathcal{N}$, in contradiction with Lemma 4.

From Lemma 3 it follows that $p_1(\delta)$ is strictly monotone in $\mathcal{D}$, and as a consequence
one can parametrize $\mathcal{N}$ through $\theta = p_1(\delta)$.
And from Remark 1 one knows that
$\inf_{\delta \in \mathcal{D}} p_1(\delta) = 0$,
hence the parameter space is $(0,p_{10}]$, where $p_{10}=\max_{\delta \in \mathcal{D}} p_1(\delta)$.
Equivalently, one can parametrize $\mathcal{N}$ through $\eta = - \log p_1(\delta)$
with parameter space $[\eta_0,\infty)$, where $\eta_0 = -\log p_{10}$.
\hfill{$\square$}

\section{Bibliography}

\begin{description}

\item[] Arnold, B.C., Balakrishnan, N., Nagaraja, H.N. (1992). \emph{A First Course in Order Statistics}. New York, Wiley.

\item[]
Consul, P.C. (1984). On the distributions of order statistics for a random sample size. \emph{Statistica Neerlandica}, 38, 249-256.

\item[]
Cordeiro, G.M., Castro, M. (2009). A new family of generalized distributions. \emph{Journal of Statistical Computation \& Simulation}, 81, 883-898.

\item[]
Cordeiro, G.M., Ortega, E.M.M., Cunha, D.C.C. (2013). The exponential generalized class of distributions. \emph{Journal of Data Science}, 11, 1-27.

\item[]
Engen, (1974). On species frequency models. \emph{Biometrika}, 61, 263-270.

\item[]
Fama, E.F., Roll, R. (1968). Some properties of symmetric stable distributions. \emph{Journal of the American Statistical Association}, 63, 817-836.

\item[]
Gupta, D., Gupta, R.C. (1984). On the distributions of order statistics for a random sample size. \emph{Statistica Neerlandica}, 38, 13-19.

\item[]
Louzada, F., Beret, E.M.P, Franco, M. A. P. (2012). On the distribution of the minimum or maximum of a
random number of i.i.d. lifetime random variables. \emph{Applied Mathematics}, 3, 350-353.

\item[]
Marshall, A.W., Olkin, I. (1997). A new method for adding a parameter to a family of distributions with application to the exponential
and Weibull families. \emph{Biometrika}, 84, 641-652.

\item[]
Rachev, S.T., Resnick, S. (1991). Max-geometric infinite divisibility and stability. \emph{Communications in Statistics. Stochastic Models,} 7, 191-218.

\item[] Raghunandanan, K. and Patil, S.A. (1972). On order statistics for random sample size. \emph{Statistica Neerlandica}, 26, 121-126.

\item[] Rohatgi, V.K. (1987). Distribution of order statistics with random sample size. \emph{Communications in Statistics.
Theory and Methods}, 16, 3739-3743.

\item[]
Shaked, M. (1975). On the distribution of the minimum and of the maximum of a random number of i.i.d. random variables.
In \emph{Statistical Distributions in Scientific Work}. Vol. I. ed. G.P. Patil, S. Kotz and J.K. Ord. Reidel, Dordrecht. pp. 363-380.

\item[]
Shaked, M., Wong, T. (1997). Stochastic comparisons of random minima and maxima. \emph{Journal of Applied Probability}, 34, 420-425

\end{description}

\end{document}